\documentclass[11pt,a4paper]{amsart}
\usepackage[dvips]{epsfig}
\usepackage{graphics}
\usepackage{latexsym}
\usepackage{verbatim}
\usepackage{amsmath}
\usepackage{amsthm}
\usepackage{amssymb}
\usepackage{MnSymbol}
\usepackage{stmaryrd}
\usepackage{hyperref}
\usepackage[]{hyperref}
\hypersetup{
    colorlinks=true,%
    filecolor=black,%
    linkcolor=black,%
    urlcolor=black
}
\usepackage [table]{xcolor}
\usepackage{multirow}
\usepackage{float}
\usepackage{tikz}
\usepackage{subfig}

\graphicspath{{pix/}}	
\newtheorem{theorem}{Theorem}[section]

\newtheorem{remark}[theorem]{Remark}

\newcommand{\RR}{\mathbb{R}}

\newcommand{\ZZ}{\mathbb{Z}}

\newcommand{\G}{\mathcal{G}}

\usepackage{algorithmic}
\newcommand{\ens}[1]{\mathbb{#1}}
\definecolor{SkyBlue}{RGB}{26 100 191}
\definecolor{GreenYellow}{RGB}{154 205 50}

\begin{document}

\title[Inverse Lax-Wendroff method for boundary conditions of Boltzmann equations]{Inverse Lax-Wendroff method for boundary conditions of Boltzmann type models}\thanks{The authors are  partially supported by the European Research Council ERC Starting Grant 2009,  project 239983-\textit{NuSiKiMo}}

\author{Francis Filbet and Chang Yang}

\hyphenation{bounda-ry rea-so-na-ble be-ha-vior pro-per-ties
cha-rac-te-ris-tic}

\maketitle

\begin{abstract}
In this paper we present a new algorithm based on a Cartesian mesh for the  numerical approximation of  kinetic models on complex geometry boundary. Due to the high dimensional property, numerical algorithms based on unstructured meshes for a complex geometry are not appropriate. Here we propose to develop an inverse Lax-Wendroff procedure, which was recently introduced  for conservation laws \cite{bibTS}, to the kinetic equations. Applications in $1D\times 3D$ and $2D\times 3D$ of this algorithm for Boltzmann type operators (BGK, ES-BGK models) are then presented and numerical results illustrate the accuracy properties of this algorithm.
\end{abstract}

\vspace{0.1cm}

\noindent 
{\small\sc Keywords.}  {\small Inverse Lax-Wendroff procedure, WENO, Boltzmann type models}

\tableofcontents

\section{Introduction} 
\setcounter{equation}{0}
\label{sec:Intro}

We are interested in the numerical approximation of  solutions to kinetic equations set in a complex geometry with different type of boundary conditions.  Unfortunately, classical structured or unstructured meshes are not  appropriate due to the high dimensional property of kinetic problem. In contrast, the Cartesian mesh makes the numerical method efficient and easy to implement. The difficulty is that obviously grid points are usually not located on the physical boundary when using a Cartesian mesh, thus a suitable numerical method to capture the boundary condition on the complex geometry is required. Several numerical methods based on Cartesian mesh have been developed in computational fluid dynamics in last decade. Among these methods, the immersed boundary method (IBM), first  introduced by Peskin~\cite{bibP} for the study of biological fluid mechanics problems, has attracted considerable attention because of its use of regular Cartesian grid and great simplification of tedious grid generation task. The basic idea of immersed boundary method is that the effect of the immersed boundary on the surrounding fluid is represented through the introduction of forcing terms in the momentum equations. In conservation laws, two major classes immersed boundary like methods  can be distinguished on different discretization types. The first class is Cartesian cut-cell method~\cite{bibICM}, which is based on a finite volume method. This conceptually simple approach ``cuts'' solid bodies out of a background Cartesian mesh. Thus we have several polygons (cut-cells) along the boundary. Then the numerical flux at the boundary of these cut-cells are imposed by using the real boundary conditions. This method satisfies well the conservation laws, however to determine the polygons is still a delicate issue. The second class is based on finite difference method. To achieve a high order interior scheme, several ghost points behind the boundary are added. For instance for solving hyperbolic conservations laws, an inverse Lax-Wendroff type procedure is used to impose some artificial values on the ghost points~\cite{bibTS}. The idea of the inverse Lax-Wendroff  procedure (ILW) is to use successively the partial differential equation  to write the normal derivatives at the inflow boundary in terms of the tangential and time derivatives of given boundary conditions. From these normal derivatives, we can obtain accurate values of ghost points using a Taylor expansion of the solution from a point located on the boundary. 

The goal of this paper is to extend the inverse Lax-Wendroff procedure to kinetic equations together with an efficient time discretization technique \cite{bibF,bibFS2010} for problems where boundary conditions play a significant role in the long time asymptotic behavior of the solution. In particular, for low speed and low Knudsen flows for which DSMC methods are unsuitable. Therefore, the  main issue relies on that the inflow is no longer a given function, while it is determined by the outflow. For this, we proceed in three steps: we first compute the outflow at the  ghost points. To maintain high order accuracy  and to prevent oscillations caused by shocks, we use a weighted essentially non-oscillatory (WENO) type extrapolation to approximate the ghost points by using the values of  interior mesh points. In the same time, we can extrapolate the outflow located at the  boundary associated with ghost points. Then, we compute the inflow at the boundary by using the outflow obtained in the first step and Maxwell's boundary conditions. Finally, we perform the inverse Lax-Wendroff procedure to approximate the inflow on the ghost points, where the key point is to replace the normal derivatives by a reformulation of the original kinetic equation.

For simplicity, we only consider simple collision operator as we adapt the ellipsoidal statistics BGK or ES-BGK model introduced by Holway~\cite{bibH}. This model gives the correct transport coefficients for Navier-Stokes approximation, so that Boltzmann or ES-BGK simulations are expected to give the same results for dense gases. Let us emphasize that Direct Simulation Monte-Carlo methods (DSMC) have been performed to the ES-BGK model in complex geometry. However DSMC approach is not computationally efficient for nonstationary or low Mach number flows  due to the requirement to perform large amounts of data simpling in order to reduce the statistical noise. In contrast, F. Filbet \& S. Jin recently proposed a deterministic  asymptotic preserving scheme for the ES-BGK model, where the entire equation can be solved explicitly and it can capture the macroscopic fluid dynamic limit even if the small scale determined by the Knudsen number is not numerically  resolved \cite{bibFS2011}. We will use this scheme to solve ES-BGK model while on the boundary the inverse Lax-Wendroff procedure will be applied.

The outline of the paper is as follows. In Section~\ref{sec:ILW} we describe precisely the inverse Lax-Wendroff procedure to Maxwell's boundary condition in 1D and 2D space dimension. Then in Section~\ref{sec:ESBGK} we present the ES-BGK model and the application of inverse Lax-Wendroff procedure to this model. In Section~\ref{sec:Num} a various numerical examples are provided  in $1D\times 3D$ and $2D\times 3D$ to demonstrate the interest and the efficiency of our method in term of accuracy and complexity. Finally a conclusion and some perspectives are given in Section~\ref{sec:conc}.

\section{Numerical method to Maxwell's boundary conditions}
\label{sec:ILW}
\setcounter{equation}{0}
The fundamental kinetic equation for rarefied gas is the Boltzmann equation
\begin{equation}
  \frac{\partial f}{\partial t}+\mathbf{v}\cdot\nabla_{\mathbf{x}}f=\frac{1}{\varepsilon}\mathcal{Q}(f),
\label{eq:boltzmann}
\end{equation}
which governs the evolution of the density $f(t,\mathbf{x},\mathbf{v})$  of monoatomic particles in the phase, where $\mathbf{x}\,\,\in\Omega_{\mathbf{x}}\subset\mathbb{R}^{d_{_\mathbf{x}}},\,\,\,\mathbf{v}\,\in\,\mathbb{R}^{3}$. The collision operator is either given by the full Boltzmann operator  

\begin{equation} \label{eq:Q2}
\mathcal{Q}(f)(\mathbf{v}) = \int_{\RR^{3}}
 \int_{\ens{S}^{2}}  B(|\mathbf{v}-\mathbf{v}_\star|, \cos \theta) \,
 \left( f^\prime_\star  f^\prime \,-\, f_\star f \right) \, d\sigma \, d\mathbf{v}_\star
 \end{equation}
or by a simplified model as the BGK or ES-BKG operator (see the next section). Boltzmann's  type collision operator has the fundamental properties of conserving mass, momentum and energy: at the formal level
\begin{equation*}
  \int_{\mathbb{R}^{3}}\mathcal{Q}(f)\;\phi(\mathbf{v})\;d\mathbf{v}=0,\,\,\,\phi(\mathbf{v})=1,\,\mathbf{v},\,{|\mathbf{v}|^2}.
\end{equation*}
Moreover, the equilibrium is the local Maxwellian distribution namely:
\begin{equation*}
  \mathcal{M}[f](t,\mathbf{x},\mathbf{v})=\frac{\rho(t,\mathbf{x})}{(2\pi \,T(t,\mathbf{x}))^{3/2}}\exp\left(-\frac{|\mathbf{u}(t,\mathbf{x})-\mathbf{v}|^2}{2\,T(t,\mathbf{x})}\right),
\end{equation*}
where $\rho$, $\mathbf{u}$, $T$ are the density, macroscopic velocity and the temperature of the gas, defined by
\begin{equation}
  \left\{
\begin{array}{l}
  \displaystyle\rho(t,\mathbf{x})\,=\,\int_{\mathbb{R}^{3}}f(t,\mathbf{x},\mathbf{v})d\mathbf{v},\\
\,
\\
  \displaystyle\mathbf{u}(t,\mathbf{x})\,=\,\frac{1}{\rho(t,\mathbf{x})}\int_{\mathbb{R}^{3}}\mathbf{v}f(t,\mathbf{x},\mathbf{v})d\mathbf{v},\\
\,
\\
  \displaystyle T(t,\mathbf{x})\,=\,\frac{1}{3\,\rho(t,\mathbf{x})} \int_{\mathbb{R}^{3}}|\mathbf{u}(t,\mathbf{x})-\mathbf{v}|^2f(t,\mathbf{x},\mathbf{v})d\mathbf{v}.
\end{array}
\right.  \label{eq:macro}
\end{equation}

In order to define completely the mathematical problem \eqref{eq:boltzmann}, suitable boundary conditions on $\partial\Omega_{\mathbf{x}}$ should be appled. Here we consider wall type boundary conditions  introduced by Maxwell~\cite{bibM}, which is assumed that the fraction $(1-\alpha)$ of the emerging  particles has been reflected elastically at the wall, whereas the remaining fraction $\alpha$ is thermalized  and leaves the wall in a Maxwellian distribution. The parameter $\alpha$ is called accommodation coefficient~\cite{bibC}.

More precisely, for $\mathbf{x}\in\partial\Omega_{\mathbf{x}}$ the smooth boundary $\partial\Omega_{\mathbf{x}}$ is assumed to have a unit inward normal $\mathbf{n}(\mathbf{x})$ and for $\mathbf{v}\cdot\mathbf{n}(\mathbf{x})\geq0$, we assume that at the solid boundary a fraction $\alpha$ of particles is absorbed by the wall and then re-emitted with the velocities corresponding to those in a still gas at the temperature of the solid wall, while the remaining portion $(1-\alpha)$ is perfectly reflected. This is equivalent to impose for the ingoing velocities
\begin{equation}
  f(t,\mathbf{x},\mathbf{v})\,=\,(1-\alpha)\,\mathcal{R}[f(t,\mathbf{x},\mathbf{v})]\,+\,\alpha\,\mathcal{M}[f(t,\mathbf{x},\mathbf{v})],\,\,\,\mathbf{x}\in\partial\Omega_{\mathbf{x}},\,\,\,\mathbf{v}\cdot\mathbf{n}(\mathbf{x})\geq0,
\label{eq:BCingoing}
\end{equation}
with $0\leq\alpha\leq1$ and
\begin{equation}
\label{op:BC}
  \left\{
  \begin{array}{lll}
    \mathcal{R}[f(t,\mathbf{x},\mathbf{v})] &=&f(t,\mathbf{x},\mathbf{v}-2(\mathbf{v}\cdot\mathbf{n}(\mathbf{x}))\mathbf{n}(\mathbf{x})),\\[3mm]
    \mathcal{M}[f(t,\mathbf{x},\mathbf{v})]&=&\mu(t,\mathbf{x})\,f_{w}(\mathbf{v}).
  \end{array}
  \right.
\end{equation}
By denoting $T_{w}$ the temperature of the solid boundary, $f_{w}$ is given by
\begin{equation}
  f_{w}(\mathbf{v})\,\,:=\,\,\exp\left(-\frac{\mathbf{v}^2}{2\,T_{w}}\right),
\end{equation}
and the value of $\mu(t,\mathbf{x})$ is determined by mass conservation at the surface of the wall for any $t\in\mathbb{R}^+$ and $\mathbf{x}\in\partial\Omega_{\mathbf{x}}$
\begin{equation}
\label{eq:MU}
  \mu(t,\mathbf{x})\int_{\mathbf{v}\cdot\mathbf{n}(\mathbf{x})\geq0}f_{w}(\mathbf{v})\mathbf{v}\cdot\mathbf{n}(\mathbf{x})d\mathbf{v}\,=\,-\int_{\mathbf{v}\cdot\mathbf{n}(\mathbf{x})<0}f(\mathbf{v})\mathbf{v}\cdot\mathbf{n}(\mathbf{x})d\mathbf{v}.
\end{equation}
This boundary condition~\eqref{eq:BCingoing} guarantees the global conservation of mass~\cite{bibF}.

 In this paper we only apply a second order finite difference method to discretize the transport term of~\eqref{eq:boltzmann} but higher order schemes \cite{bibJS} may be applied. Then to keep the order of accuracy of the  method, two ghost points should be added in each spatial direction.   To impose  $f$ at the ghost points, we will  apply the inverse Lax-Wendroff procedure proposed in~\cite{bibTS} for conservation laws.

Suppose that   the distribution function $f$ at time level $t^n$ for all interior points are already known, we now construct $f$ at the ghost points.


\subsection{One-dimensional case in space}
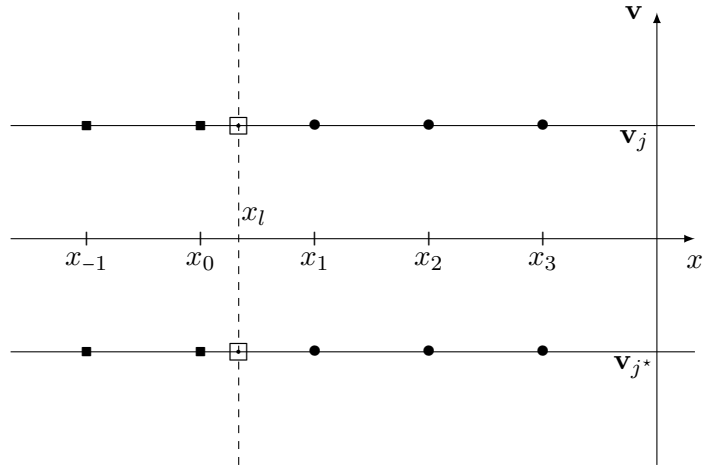
\begin{figure}[h]
  \begin{center} 
   \begin{tikzpicture}[point/.style={circle,inner sep=0pt,minimum size=2pt,fill=black}]
      \begin{scope}[>=latex]          
        \draw[->] (-5,0) --(4,0)node[midway,right] {};
        \draw[-] (-5,1.5) --(4,1.5)node[midway,right] {};
        \draw[-] (-5,-1.5) --(4,-1.5)node[midway,right] {};
        \draw[->] (3.5,-3) --(3.5,3)node[midway,right] {};
        \draw[dashed] (-2,-3) --(-2,3)node[midway,right] {};
      \end{scope}
      
       \draw   (-4,-0.3) node{$x_{-1}$};
       
       \draw   (-2.5,-0.3) node{$x_0$};
       
       \draw   (-1.8,0.3) node{$x_l$};
       
       \draw   (-1,-0.3) node{$x_{1}$};
      
       \draw   (0.5,-0.3) node{$x_{2}$};
        \draw   (2,-0.3) node{$x_{3}$};

       \draw   (-4,1.48) node{{$\filledsquare$}};
       \draw   (-2.5,1.48) node{{$\filledsquare$}};
       \draw   (-2,1.5) node{$\boxdot$};
       \draw   (-1,1.5) node{$\bullet$};
       \draw   (0.5,1.5) node{$\bullet$};
       \draw   (2,1.5) node{$\bullet$};

       \draw   (-4,0) node{$\medvert$};
       \draw   (-2.5,0) node{$\medvert$};
       
       \draw   (-1,0) node{$\medvert$};

       \draw   (0.5,0) node{$\medvert$};
        \draw   (2,0) node{$\medvert$};

       \draw   (-4,-1.51) node{{$\filledsquare$}};
       \draw   (-2.5,-1.51) node{{$\filledsquare$}};
       \draw   (-2,-1.5) node{$\boxdot$};
      
       \draw   (-1,-1.5) node{$\bullet$};
      
       \draw   (0.5,-1.5) node{$\bullet$};
       \draw   (2,-1.5) node{$\bullet$};
      
       \draw   (4,-0.3) node{$x$};
       \draw   (3.2,3) node {$\mathbf{v}$};
       \draw   (3.2,1.3) node {$\mathbf{v}_{j}$};
       \draw   (3.2,-1.7) node {$\mathbf{v}_{j^\star}$};
    \end{tikzpicture}
\caption{\label{fig:1Ddomain}A portion of mesh in spatially one dimensional case. $\bullet$ is interior point, $\filledsquare$ is ghost point, $\boxdot$ is the left hand side boundary.}
  \end{center}
\end{figure}

We start with spatially one-dimensional  problem, that is $d_{\mathbf{x}}=1$. In this case the Boltzmann equation reads:
\begin{equation}
    \begin{array}{ll}
      \displaystyle\frac{\partial f}{\partial t}+v_x\frac{\partial f}{\partial x}=\frac{1}{\varepsilon}\mathcal{Q}(f),&(x,\mathbf{v})\in [x_{l},x_{r}]\times\mathbb{R}^{3},
    \end{array}
\label{eq:1D}
\end{equation}
where $x_l$ and $x_r$ are the left and right boundaries respectively, $v_x$ is the component of phase space corresponding to $x$-direction. For  the boundary condition in spatially one-dimensional  case,  the inward normal on the boundary in~\eqref{eq:BCingoing} is
\begin{eqnarray*}
  \mathbf{n}(x_l)\,=\,\begin{pmatrix}1\\0\\0 \end{pmatrix},\,\,\,
  \mathbf{n}(x_r)\,=\,\begin{pmatrix}-1\\0\\0 \end{pmatrix}.
\end{eqnarray*}
To  implement the numerical method, we assume the computational domain is a limited domain $[x_{\min},x_{\max}]\times[-V,V]^{3}$, where $(x_l,x_r)\subset[x_{\min},x_{\max}]$. The computational domain is covered by a uniform Cartesian mesh $\mathbf{X}_h\times \mathbf{V}_h$,
\begin{equation}
\left\{
\begin{array}{l}
  \mathbf{X}_h=\left\{x_{\min}=x_0\leq\cdots\leq x_i\leq\cdots\leq x_{n_x}=x_{\max}\right\},
\\
\,
\\
 \mathbf{V}_h=\left\{\mathbf{v}_j\,=\,j\,\Delta v, \quad  j=(j_1,\ldots,j_{3}) \in \ZZ^{3}, \quad |j|\leq n_v \right\}.
\end{array}
\right.
\label{eq:1Dmesh}
\end{equation}
with the mesh size $\Delta x$ and $\Delta v$ for space and velocity respectively. We only consider numerical method of ghost points near the left hand side boundary, since the  procedure for right hand side boundary is the same. Figure~\ref{fig:1Ddomain} illustrates  a portion of mesh near left boundary $x_l$, which is located between $x_0$ and $x_1$.

We construct $f$ at each ghost point in following three steps: we perform an extrapolation of $f$ to compute a high order approximation of the outflow. Then, we compute an approximation of the distribution function at the boundary using Maxwell's boundary conditions. Finally, we apply the inverse Lax-Wendroff procedure for the inflow.

 \subsubsection{First step: Extrapolation of $f$ for the outflow} At time $t=t^n$ we consider the outflow near the point $x_l$, that is  $f(t,x_l,\mathbf{v}_{j})$ where $v_{j_1}<0$.  We denote by $f_{i,j}$ an approximation of $f$ at $(x_i,\mathbf{v}_{j})$

A natural idea is to extrapolate $f$ at the left boundary $x_l$  or the ghost points $x_{0}$ and $x_{-1}$ using the values of $f$ on interior  points. For example  from the values $f_{1,j}$, $f_{2,j}$ and $f_{3,j}$, we can construct a  Lagrange polynomial  $p_2(x)\in\mathbb{P}_2(\mathbb{R})$.  Then by injecting $x_l$, $x_0$ or $x_{-1}$ into $p_2(x)$, we obtain the approximations of $f$ at the ghost points and left boundary, {\it i.e.}  $f_{l,j}$, $f_{0,j}$ and $f_{-1,j}$. However, when a shock goes out of the boundary, the high order extrapolation may lead to a severe oscillation near the shock. To prevent this, we would like have a lower order accurate but more robust extrapolation. Therefore, a WENO type extrapolation~\cite{bibTS} will be applied  and described below (see subsection~\ref{sec:WENO}) for this purpose.

 \subsubsection{Second step: Compute boundary conditions at the boundary}
In the previous step, the outflow at the boundary is obtained by extrapolation. To compute the vlaues of $f$ at the inflow boundary, we apply the Maxwell's boundary condition~\eqref{eq:BCingoing}, {i.e.}
\begin{equation}
  f_{l,j}\,=\,(1-\alpha)\;\mathcal{R}[f_{l,j}]\,+\,\alpha\,\mathcal{M}[f_{l,j}].
\end{equation}
On the one hand the specular reflection portion is given straightly by the outflow at the left boundary, which is
\begin{equation*}
  \mathcal{R}[f_{l,j}]\,=\, f_{l,j^\star}, \quad{\rm where}\quad j^\star=(-j_1,j_2,\ldots,j_{3}). 
\end{equation*}
On the other hand the diffuse one is computed by a half Maxwellian
\begin{equation*}
  \mathcal{M}[f_{l,j}]\,\,=\,\,\mu_l\,\exp\left(-\frac{|\mathbf{v}_{j}|^2}{2T_l}\right),
\end{equation*}
where $T_l$ is the given temperature at the left wall and $\mu_l$ is given by
\begin{equation*}
  \mu_l\, {\sum_{\mathbf{v}_j\cdot \mathbf{n}(x_l) \geq 0} \mathbf{v}_j\cdot \mathbf{n}(x_l) \, \exp\left(-\frac{|\mathbf{v}_j|^2}{2\,T_l}\right) } \,\,=\,\, {-\sum_{\mathbf{v}_j\cdot \mathbf{n}(x_l) \leq 0} \mathbf{v}_j\cdot \mathbf{n}(x_l) \, f_{l,j}}.
\end{equation*}

 \subsubsection{Third step: Approximation of $f$ at the inflow boundary} 
Finally we compute the values of $f$ at the ghost points for the inflow boundary. Here we cannot approximate  $f$ by an extrapolation, since the distribution function at interior points  cannot predict the inflow. Then we  extend the inverse Lax-Wendroff type procedure recently proposed in \cite{bibHSZ, bibTS, bibXZZS} for solving kinetic equations. At the left boundary $x_l$, a first order Taylor expansion gives
\begin{equation*}
  f_j(x)=f_{l,j}+(x-x_l)\left.\frac{\partial f}{\partial x}\right|_{x=x_l}+O(\Delta x^2).
\end{equation*}
Hence a first order approximation of $f$ at ghost points is
\begin{equation}
  f_{s,j}\,\,=\,\,f_{l,j}\,+\,(x_s-x_l)\,\left.\frac{\partial f}{\partial x}\right|_{x=x_l},\quad s=-1,0.
\end{equation}
We already have $f_{l,j}$ in the second step, thus it remains to obtain an approximation of the first derivative. By reformulating~\eqref{eq:1D}, we have
\begin{equation}
  \left.\frac{\partial f}{\partial x}\right|_{x=x_l}=\left.\frac{1}{v_x}\left(-\frac{\partial f}{\partial t}+\frac{1}{\varepsilon}\mathcal{Q}(f)  \right)\right|_{x=x_l}.
\label{eq:1Dreformulation}
\end{equation}
Now instead of approximating the first derivative $\partial_x f|_{x=x_l}$, we compute the time derivative $\partial_t f|_{x=x_l}$ and the collision operator $\mathcal{Q}(f)|_{x=x_l}$. An approximation of the time derivative can be computed by using several $f_{i,j}$   at previous time levels. Different approximation are obtained either a first order approximation reads
\begin{equation*}
  \left.\frac{\partial f}{\partial t}\right|_{x=x_l}=\frac{f^n_{l,j}-f^{n-1}_{l,j}}{\Delta t},
\end{equation*}
where $\Delta t$ is the time step, or one can use a WENO type extrapolation to approximate the time derivative (see subsection~\ref{sec:WENO} below). 

The last term $\mathcal{Q}(f)|_{x=x_l}$ can be computed explicitly by using $f_{l,j}$ obtained in previous two steps. Clearly this procedure is independent of the values of $f$ at interior points.

\begin{remark}
Let us observe that when $\alpha=0$ we have a pure specular reflection boundary condition. A mirror procedure can be used to approximate $f$ at the ghost points. More precisely, by considering the boundary as a mirror, we approximate the distribution at the ghost points $f(x_s,v_j)$ as
\begin{equation*}
  f(x_s,\mathbf{v}_j)\,=\,f(2\,x_l-x_s,\mathbf{v}_{j^\star}), \quad{\rm where}\quad j^\star = (-j_1,j_2,\ldots,j_{3})
\end{equation*}
where $2\,x_l\,-\,x_s$ is the mirror image point of $x_s$. Since $2\, x_l-x_s$ is located in  interior domain, we can approximate $f(2\,x_l-x_s,\mathbf{v}_{j^\star})$ by an interpolation procedure.
\end{remark}

\subsection{Two-dimensional case in space}
\label{subsec:2D}
The previous approach  can be generalized to spatially two-dimensional problem. We assume  $d_{\mathbf{x}}=2$ in equation~\eqref{eq:boltzmann}
\begin{equation}
      \displaystyle\frac{\partial f}{\partial t}+v_x\frac{\partial f}{\partial x}+v_y\frac{\partial f}{\partial y}=\frac{1}{\varepsilon}\mathcal{Q}(f),
\label{eq:2D}
\end{equation}
where the distribution function $f(t,\mathbf{x},\mathbf{v})$ is defined in $(t,\mathbf{x},\mathbf{v})\in\mathbb{R}^+\times\Omega\times\mathbb{R}^{3}$ with $\mathbf{x}=(x,y)$. We consider a computational domain $[x_{\min},x_{\max}]\times[y_{\min},y_{\max}]\times[-V,V]^{3}$,  such that $\Omega\subset[x_{\min},x_{\max}]\times[y_{\min},y_{\max}]$ and  $f(t,\mathbf{x},\mathbf{v})\approxeq 0$, for all $\|\mathbf{v}\|\geq V$. 

\begin{figure}[h]
  \begin{center} 
   \begin{tikzpicture}[point/.style={circle,inner sep=0pt,minimum size=2pt,fill=black}]
      \begin{scope}[>=latex]          
        \draw[->] (-1,0) --(7,0)node[right] {$x$};
        \draw[-] (-1,1.5) --(7,1.5)node[midway,right] {};
        \draw[-] (-1,3) --(7,3)node[midway,right] {};
        \draw[-] (-1,4.5) --(7,4.5)node[midway,right] {};
        \draw[-] (-1,6) --(7,6)node[midway,right] {};
        \draw[->] (0,-1) --(0,7)node[left] {$y$};
        \draw[-] (1.5,-1) --(1.5,7)node[midway,right] {};
        \draw[-] (3,-1) --(3,7)node[midway,right] {};
        \draw[-] (4.5,-1) --(4.5,7)node[midway,right] {};
        \draw[-] (6,-1) --(6,7)node[midway,right] {};
        \draw[dashed] (7,-0.3) --(-1,4)node[midway,right] {};
        \draw[->] (1.2,1) --(4.8,7)node[right] {$\mathbf{n}$};
        \draw[dotted] (-0.33,1) --(2,4.5)node[right] {};
        \draw[dotted] (1.3,-0.3) --(4,4)node[right] {};
        \draw[->] (3.5,4.8) .. controls +(0.2,0.1) and +(-0.1,0.2)  .. (4.,4.5)node[midway, above] {$\theta$};
      \end{scope}
      \draw   (1.5,0) node{$\filledsquare$};
      \draw   (3,0) node{$\filledsquare$};
      \draw   (4.5,0) node{$\filledsquare$};
      \draw   (6,0) node{$\filledsquare$};
      \draw   (6,6) node{$\largecircle$};
      \draw   (3,6) node{$\largecircle$};
      \draw   (4.5,6) node{$\largecircle$};
      \draw   (0,1.5) node{$\filledsquare$};
      \draw   (1.5,1.5) node{{\color{red}$\filledsquare$}};
      \draw   (3,1.5) node{$\filledsquare$};
      \draw   (4.5,1.5) node{$\bullet$};
      \draw   (6,1.5) node{$\bullet$};
      \draw   (1.5,4.5) node{$\largecircle$};
      \draw   (3,4.5) node{$\largecircle$};
      \draw   (4.5,4.5) node{$\largecircle$};
      \draw   (1.8,1.2) node{$\mathbf{x}_g$};
      \draw   (0,3) node{$\filledsquare$};
      \draw   (1.5,3) node{$\bullet$};
      \draw   (3,3) node{$\bullet$};
      \draw   (4.5,3) node{$\bullet$};
      \draw   (6,3) node{$\bullet$};
      \draw   (1.5,3) node{$\largecircle$};
      \draw   (3,3) node{$\largecircle$};
      \draw   (4.5,3) node{$\largecircle$};
      \draw   (2.1,2.05) node{$\mathbf{x}_p$};
      \draw   (0,4.5) node{$\bullet$};
      \draw   (1.5,4.5) node{$\bullet$};
      \draw   (3,4.5) node{$\bullet$};
      \draw   (4.5,4.5) node{$\bullet$};
      \draw   (6,4.5) node{$\bullet$};
     
      \draw   (0,6) node{$\bullet$};
      \draw   (1.5,6) node{$\bullet$};
      \draw   (3,6) node{$\bullet$};
      \draw   (4.5,6) node{$\bullet$};
      \draw   (6,6) node{$\bullet$};

      \draw   (2,2.4) node{$\boxdot$};
       
      \draw   (0,-1.3) node{$i_x-2$};
      \draw   (1.5,-1.3) node{$i_x-1$};
      \draw   (3,-1.3) node{$i_x$};
      \draw   (4.5,-1.3) node{$i_x+1$};
      \draw   (6,-1.3) node{$i_x+2$};
      \draw   (-1.5,0) node{$i_y-2$};
      \draw   (-1.5,1.5) node{$i_y-1$};
      \draw   (-1.5,3) node{$i_y$};
      \draw   (-1.5,4.5) node{$i_y+1$};
      \draw   (-1.5,6) node{$i_y+2$};
      \draw   (0.95,2.95) node{$\boxdot$};
      \draw   (2.7,2) node{$\boxdot$};
    \end{tikzpicture}  
\caption{\label{fig:2Ddomain}Spatially two-dimensional Cartesian mesh. $\bullet$ is interior point, $\filledsquare$ is ghost point, $\boxdot$ is the point at the boundary, $\largecircle$ is the point for extrapolation, the dashed line is the boundary.}
  \end{center}
\end{figure}
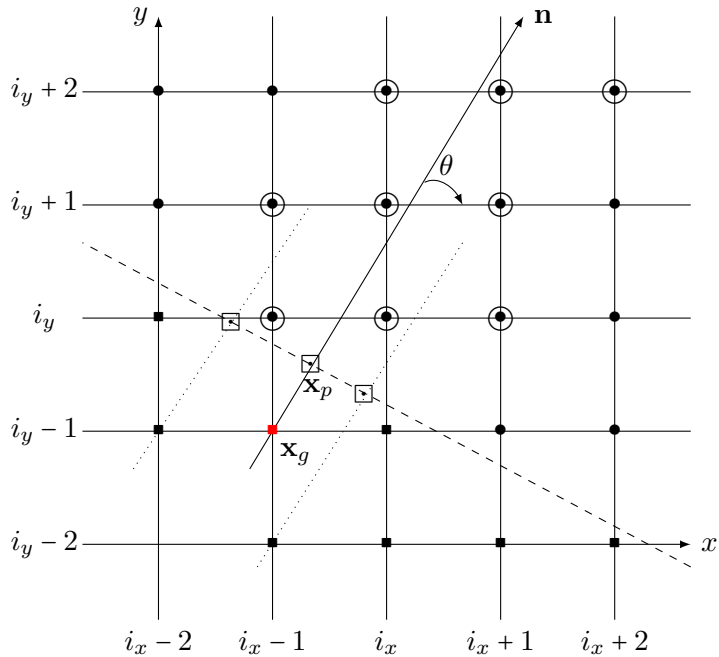

The computational domain is covered by an uniform Cartesian mesh $\mathbf{X}_h\times \mathbf{V}_{h}$, where $\mathbf{X}_h$, $\mathbf{ V}_{h}$ are defined similarly to \eqref{eq:1Dmesh}. The mesh steps are respectively $\Delta x$, $\Delta y$ and $\Delta v$. In Figure~\ref{fig:2Ddomain}, we present a portion of spatial mesh near the boundary. From a ghost point $\mathbf{x}_g$, we can find an inward normal $\mathbf{n}$, which  crosses the boundary at $\mathbf{x}_p$.


For the 2D case in space, the numerical approximation of the distribution function $f$ at ghost points is similar to the one dimensional case. However, there are two major differences. First to compute $\mathcal{R}[f]$ in the second step, the corresponding outflow may not locate on phase space mesh. Secondly to approximate the normal derivative in the third step, besides the time derivative and collision operator we need also the tangential derivative at $\mathbf{x}_p$.  Once again, we present the method in three steps:

\subsubsection{First step: Extrapolation of $f$ for outflow}
Let us assume that the values of the distribution function $f$ on the grid points in $\Omega$ are given. To approximate $f$ at a ghost point, for instance $\mathbf{x}_s$, we first construct a stencil $\mathcal{E}$ composed of grid points of $\Omega$ for the extrapolation. For instance as it is shown in Figure~\ref{fig:2Ddomain},  the inward normal $\mathbf{n}$ intersects the grid lines $y=y_{i_y}$, $y_{i_y+1}$, $y_{i_y+2}$. Then we choose the three nearest point of the cross point in each line, {\it i.e.} marked by a large circle. From these nine points, we can build a  Lagrange polynomial $q_2(\mathbf{x})\in\mathbb{Q}_2(\mathbb{R}^2)$. Therefore  we evaluate the polynomial $q_2(\mathbf{x})$ at $\mathbf{x}_s$ or $\mathbf{x}_p$, and obtain an approximation of $f$ at the boundary and at ghost points. As for the 1D case, a WENO type extrapolation can be used to prevent spurious oscillations, which will be detailed in subsection~\ref{sec:WENO}.


\subsubsection{Second step: Compute boundary conditions at the boundary}
In the previous step, we have obtained the outflow $f(\mathbf{x}_p, \mathbf{v}\cdot\mathbf{n}<0)$ at the boundary $\mathbf{x}_p$. By using~\eqref{eq:BCingoing} as we did for the 1D case, we can similarly compute the distribution  function $f$ for $\mathbf{v}\cdot\mathbf{n}\geq 0$. However this time to compute the distribution function for specular reflection
$$
\mathcal{R}[f(\mathbf{x}_p,\mathbf{v})]\,=\, f(\mathbf{x}_p,\mathbf{v}-2(\mathbf{v}\cdot\mathbf{n})\mathbf{n}),\quad \forall\mathbf{v}\in\mathbf{V}_{h},
$$ 
the vector fields $\mathbf{v} - 2\,(\mathbf{v}\cdot\mathbf{n})\,\mathbf{n}$ may not be located on a grid point. Therefore, we  interpolate $f$ in phase space $(\mathbf{x}_p,\mathbf{v}-2(\mathbf{v}\cdot\mathbf{n})\mathbf{n})$ using the values computed from the outflow $f(\mathbf{x}_p,\mathbf{v})$ such that $\mathbf{v}\cdot\mathbf{n}\geq 0$.

\subsubsection{Third step: Approximation of $f$ at the inflow boundary}
We have obtained the values of  $f$ at the boundary points $\mathbf{x}_p$ for all $\mathbf{v}\in\mathbf{V}_{h}$ in previous two steps. Now we reconstruct the values of $f$  for the velocity grid points such that $\mathbf{v}\cdot\mathbf{n}\leq 0$ at the ghost point $\mathbf{x}_g$ by a simple Taylor expansion in the inward normal direction. To this end, we set up a local coordinate system at $\mathbf{x}_p$ by
\begin{eqnarray*}
  \hat{\mathbf{x}}=\begin{pmatrix}
    \hat x\\\hat y
  \end{pmatrix}
=
  \begin{pmatrix}
    \cos\theta&\sin\theta\\
    -\sin\theta&\cos\theta
  \end{pmatrix}
  \begin{pmatrix}
    x\\y
  \end{pmatrix},
\end{eqnarray*}
where $\theta$ is the angle between the inward normal $\mathbf{n}$ and the $x$-axis  illustrated in Figure~\ref{fig:2Ddomain}. Thus the first order approximation of $f(\mathbf{x}_g,\mathbf{v})$ reads
\begin{equation*}
  f(\mathbf{x}_g,\mathbf{v})\approxeq \hat f(\hat {\mathbf{x}}_p,\mathbf{v})+(\hat x_g-\hat x_p)\frac{\partial \hat f}{\partial\hat x}(\hat{\mathbf{x}}_p,\mathbf{v}),
\end{equation*}
where $\hat f(\hat{\mathbf{x}}_p,\mathbf{v})=f(\mathbf{x}_p,\mathbf{v})$ and $\frac{\partial \hat f}{\partial\hat x}(\hat{\mathbf{x}}_p,\mathbf{v})$ is the first order normal derivative at the boundary $\mathbf{x}_p$. To approximate $\frac{\partial \hat f}{\partial\hat x}(\hat{\mathbf{x}}_p,\mathbf{v})$, we use inverse Lax-Wendroff procedure. Firstly, we rewrite the equation~\eqref{eq:2D} in the local coordinate system as
\begin{equation}
  \displaystyle\frac{\partial \hat f}{\partial t}+\hat v_x\frac{\partial \hat f}{\partial \hat  x}+\hat v_y\frac{\partial \hat f}{\partial \hat  y}=\frac{1}{\varepsilon}\mathcal{Q}(\hat f),
  \label{eq:2D:changevari}
\end{equation}
where $\hat v_x=v_x\cos\theta+v_y\sin\theta$, $\hat v_y=-v_x\sin\theta+v_y\cos\theta$. Then a reformulation of~\eqref{eq:2D:changevari} yields
\begin{equation}
  \frac{\partial \hat f}{\partial\hat x}(\hat{\mathbf{x}}_p,\mathbf{v})=-\frac{1}{\hat v_x}\left.\left(\frac{\partial \hat f}{\partial t}+\hat v_y\frac{\partial \hat f}{\partial \hat y}-\frac{1}{\varepsilon}\mathcal{Q}(\hat f)\right)\right|_{\hat{\mathbf{x}}=\hat{\mathbf{x}}_p}.
  \label{eq:2Dreformulation}
\end{equation}
Finally instead of approximating $\frac{\partial \hat f}{\partial\hat x}(\hat{\mathbf{x}}_p,\mathbf{v})$ directly, we approximate the time derivative $\frac{\partial \hat f}{\partial t}$, tangential derivative $\frac{\partial \hat f}{\partial \hat y}$ and  collision operator $\mathcal{Q}(\hat f)$. Similarly as in spatially 1D case, we  compute $\frac{\partial \hat f}{\partial t}$ and  $\mathcal{Q}(\hat f)$. It remains to approximate  $\frac{\partial \hat f}{\partial \hat y}$. For this, some neighbor points of $\mathbf{x}_p$ at the boundary  are required (See the empty squares in Figure~\ref{fig:2Ddomain}).  We then perform an essentially non-oscillatory (ENO) procedure~\cite{bibHEOC}  for this numerical differentiation to avoid the discontinuity.

\subsection{WENO type extrapolation}
\label{sec:WENO}
A WENO type extrapolation~\cite{bibTS} was developed to prevent oscillations and maintain accuracy. The key point of WENO type extrapolation is to define smoothness indicators, which is designed to help us  choose automatically between the high order accuracy  and the low order but more robust extrapolation. Here we describe this method in spatially 1D and 2D cases. Moreover we will give a slightly modified version of the  method such that  the smoothness indicators are invariant with respect to the scaling of $f$.

\subsubsection{One-dimensional WENO type extrapolation} Assume that we have a stencil of three points $\mathcal{E}=\{x_1,x_2,x_3\}$ showed in Figure~\ref{fig:1Ddomain} and denote the corresponding distribution function by $f_1$, $f_2$, $f_3$. Instead of extrapolating $f$ at ghost point $x_g$ by Lagrange polynomial, we use following Taylor expansion
\begin{equation*}
  f_g=\sum_{k=0}^2\frac{(x_g-x_l)^2}{k!}\left.\frac{d^kf}{dx^k}\right|_{x=x_l}
\end{equation*}
 We aim to obtain a $(3-k)$-th order approximation of $\left.\frac{d^kf}{dx^k}\right|_{x=x_l}$ denoted by $f_l^{(k)}$, $k=0,1,2$. Three candidate substencils are given by
\begin{equation*}
  S_r=\{x_1,\dots,x_{r+1}\},\,\,r=0,1,2.
\end{equation*}
In each substencil $S_r$, we could construct a  Lagrange polynomial  $p_r(x)\in\mathbb{P}_r(\mathbb{R})$
\begin{equation*}
\left\{
\begin{array}{l}
  \displaystyle p_0(x)\,=\,f_1,\\\,\\
  \displaystyle p_1(x)\,=\, f_1+\frac{f_{2}-f_1}{\Delta x}(x-x_1),\\\,\\
 \displaystyle p_2(x) \,=\, f_1+\frac{f_{2}-f_1}{\Delta x}(x-x_1)+\frac{f_3-2f_2+f_1}{2\Delta x^2}(x-x_1)(x-x_2).
\end{array}\right.
\end{equation*}
 We now look for the WENO type extrapolation in the form
\begin{equation*}
  f^{(k)}_l=\sum_{r=0}^2w_r\frac{d^kp_r(x)}{dx^k}(x_l),
\end{equation*}
where $w_r$ are the nonlinear weights depending on $f_i$. We expect  that  $f^{(k)}_l$ has $(3-k)$-order accurate in the case $f(x)$ is smooth in $S_2$. The nonlinear weights are given by
\begin{equation*}
  w_r=\frac{\alpha_r}{\sum^2_{s=0}\alpha_s},
\end{equation*}
with 
\begin{equation*}
  \alpha_r=\frac{d_r}{(\varepsilon+\beta_s)^2},
\end{equation*}
where $\varepsilon=10^{-6}$ and $\beta_r$ are the new smoothness indicators determined by
\begin{eqnarray*}
  \beta_0&=&\Delta x^2,\\ \,\\
  \beta_1&=&\frac{1}{\varepsilon+f_1^2+f_2^2}\sum^2_{l=1}\int_{x_{0}}^{x_1}\Delta x^{2l-1}\left(\frac{d^l}{dx^l}p_1(x)\right)^2dx\\
  &=&\frac{(f_2-f_1)^2}{\varepsilon+f_1^2+f_2^2}\\ \,\\
  \beta_2&=&\frac{1}{\varepsilon+f_1^2+f_2^2+f^2_3}\sum^2_{l=1}\int_{x_{0}}^{x_1}\Delta x^{2l-1}\left(\frac{d^l}{dx^l}p_2(x)\right)^2dx\\
  &=&\frac{1}{12(\varepsilon+f_1^2+f_2^2+f_3^2)}(61f_1^2+160f_2^2+25f_3^2+74f_1f_3-192f_1f_2-124f_2f_3).
\end{eqnarray*}
We remark that the smoothness indicators $\beta_1$ and $\beta_2$ have the factors $\frac{1}{\varepsilon+\sum_{m=1}^{r+1}f_m^2}$, which guarantee that the indicators are invariant of the scaling of $f_i$.

\subsubsection{Two-dimensional extrapolation}
The two-dimensional extrapolation is a straightforward expansion of 1D case. The substencils $S_r,\,\,r={0,1,2}$ for extrapolation are chosen around the inward normal $\mathbf{n}$ such that we can construct Lagrange polynomial of degree $r$. For instance in  Figure~\ref{fig:2Ddomain}, the three substencils are respectively  
$$
\left\{\begin{array}{l}
 \displaystyle S_0\,=\;\{(x_{i_x},y_{i_y})\},
\\ \,\\
 \displaystyle S_1\,=\,\{(x_{i_x-1},y_{i_y}),(x_{i_x},y_{i_y}),(x_{i_x},y_{i_y+1}),(x_{i_x+1},y_{i_y+1})\},
\\\,\\
 \displaystyle S_2\,=\,  \{(x_{i_x-1},y_{i_y}),(x_{i_x},y_{i_y}),(x_{i_x+1},y_{i_y}),(x_{i_x-1},y_{i_y+1}),
\\\,\\
 \displaystyle\,\quad (x_{i_x},y_{i_y+1}),(x_{i_x+1},y_{i_y+1}),(x_{i_x},y_{i_y+2}),(x_{i_x+1},y_{i_y+2}),(x_{i_x+2},y_{i_y+2})\}.
\end{array}\right.
$$
Once the substencils $S_r$ are chosen, we could easily construct the Lagrange polynomials in $\mathbb{Q}_r(\mathbb{R}^2)$
\begin{equation*}
  q_r(\mathbf{x})=\sum_{m=0}^r\sum_{l=0}^ra_{l,m}x^ly^m
\end{equation*}
satisfying 
\begin{equation*}
  q_r(\mathbf{x})=f(\mathbf{x}),\,\,\mathbf{x}\in S_r.
\end{equation*}
Then the  WENO extrapolation has the form
\begin{equation}
  f(\mathbf{x})=\sum_{r=0}^2w_rq_r(\mathbf{x}),\,\,\mathbf{x}\in S_r,
\end{equation}
where $w_r$ are the nonlinear weights, which are chosen to be
\begin{equation*}
  w_r=\frac{\alpha_r}{\sum_{s=0}^2\alpha_s},
\end{equation*}
with 
\begin{equation*}
  \alpha_r=\frac{d_r}{(\varepsilon+\beta_r)^2},
\end{equation*}
where $\varepsilon=10^{-6}$, $d_0=\Delta x^2+\Delta y^2$, $d_1=\sqrt{\Delta x^2+\Delta y^2}$, $d_2=1-d_0-d_1$. $\beta_r$ are the smoothness indicators determined by
\begin{eqnarray*}
  \beta_0&=&\Delta x^2+\Delta y^2,\\
  \beta_r&=&\frac{1}{\varepsilon+\sum_{\mathbf{x}\in S_r}f(\mathbf{x})^2}\sum_{1\leq|\alpha|\leq r}\int_K|K|^{|\alpha|-1}(D^\alpha q_r(\mathbf{x}))^2d\mathbf{x},\,\,r=1,2,
\end{eqnarray*}
where $\alpha$ is a multi-index and $K=[x_p-\Delta x/2,x_p+\Delta x/2]\times[y_p-\Delta y/2,y_p+\Delta y/2]$, $\mathbf{x}_p=(x_p,y_p)$.
\section{Application to the ES-BGK model}
\label{sec:ESBGK}
\setcounter{equation}{0}
The Boltzmann equation~\eqref{eq:boltzmann} governs well the evolution of density $f$ in kinetic regime and also in the continuum regime~\cite{bibFS2011}. However the quadratic collision operator $\mathcal{Q}(f)$ has a rather complex form such that it is very difficult to compute. Hence different simpler models have been introduced. The simplest model is the so-called BGK model~\cite{bibBGK}, which is mainly a relaxation towards a Maxwellian equilibrium state
\begin{equation}
  \mathcal{Q}(f)=\frac{\tau}{\varepsilon}(\mathcal{M}[f]-f),
\label{eq:BGKoperator}
\end{equation}
where $\tau$ depends on macroscopic quantities $\rho$ and $T$.

Although it describes the right hydrodynamical limit, the BGK model does not give the Navier-Stokes equation with correct transport coefficients in the Chapman-Enskog expansion. Holway {\it et al.}~\cite{bibH} proposed the ES-BGK model, where the Maxwellian $\mathcal{M}[f]$ in the relaxation term of ~\eqref{eq:BGKoperator} is replaced by an anisotropic Gaussian $\mathcal{G}[f]$.
 This model has correct conservation laws, yields the Navier-Stokes approximation via the Chapman-Enskog expansion with a Prandtl number less than one, and yet is endowed with the  entropy condition~\cite{bibATPP}. In order to introduce the Gaussian model, we need further notations. Define the opposite of the stress tensor
\begin{equation}
  \Theta(t,\mathbf{x})=\frac{1}{\rho}\int_{\mathbb{R}^{3}}(\mathbf{v}-\mathbf{u})\otimes(\mathbf{v}-\mathbf{u})f(t,\mathbf{x},\mathbf{v})d\mathbf{v}.
  \label{eq:stress}
\end{equation}
Therefore the translational temperature is related to the $T=\text{tr}(\Theta)/3$. We finally introduce the corrected tensor
\begin{equation*}
  \mathcal{T}(t,\mathbf{x})=[(1-\nu)T\text{I}+\nu\Theta](t,\mathbf{x}),
\end{equation*}
which can be viewed as a linear combination of the initial stress tensor $\Theta$ and of the isotropic stress tensor $T\text{I}$ developed by a Maxwellian distribution, where $\text{I}$ is the identity matrix.

The ES-BGK model introduces a corrected BGK collision operator by replacing the local equilibrium Maxwellian by the Gaussian $\mathcal{G}[f]$ defined by
\begin{equation*}
  \mathcal{G}[f]=\frac{\rho}{\sqrt{\text{det}(2\pi\mathcal{T})}}\exp\left(-\frac{(\mathbf{v}-\mathbf{u})\mathcal{T}^{-1}(\mathbf{v}-\mathbf{u})}{2}\right).
\end{equation*}
Thus, the corresponding collision operator is now 
\begin{equation}
  \mathcal{Q}(f)=\frac{\tau}{\varepsilon}(\mathcal{G}[f]-f),
\label{eq:ESBGKoperator}
\end{equation}
where $\tau$ depends on $\rho$ and $T$, the parameter $-1/2\leq\nu<1$ is used to modify the value of the Prandtl number through the formula
\begin{equation*}
  \frac{2}{3}\leq\text{Pr}=\frac{1}{1-\nu}\leq+\infty.
\end{equation*}
It  follows from the above definitions that 
\begin{equation}
\left\{\begin{array}{ll}
\displaystyle{\int_{\mathbb{R}^{3}} f(\mathbf{v})\,d\mathbf{v}  = \int_{\mathbb{R}^{3}} \G[f](\mathbf{v})\,d\mathbf{v}  = \rho,} &
\\
\;
\\
\displaystyle{ \int_{\mathbb{R}^{3}}\mathbf{v}\, f(\mathbf{v})\,d\mathbf{v} = \int_{\mathbb{R}^{3}} \mathbf{v}\,\G[f](\mathbf{v})\,d\mathbf{v}  = \rho\,\mathbf{u},} &
\\
\,
\\
\displaystyle{ \int_{\mathbb{R}^{3}} \frac{|\mathbf{v}|^2}{2}\,f(\mathbf{v})\,d\mathbf{v}  = \int_{\mathbb{R}^{3}} \frac{|\mathbf{v}|^2}{2}\,\G[f](\mathbf{v})\,d\mathbf{v}  = E} &
\end{array}\right.
\label{eq:ESBGKmacro}
\end{equation}
and
$$
\left\{
\begin{array}{ll}
\displaystyle \int_{\mathbb{R}^{3}} (\mathbf{v}-\mathbf{u}) \otimes (\mathbf{v}-\mathbf{u})\, f(\mathbf{v})\,d\mathbf{v} \,=\, \rho\, \Theta,&
\\
\,
\\
\displaystyle\int_{\mathbb{R}^{3}} (\mathbf{v}-\mathbf{u}) \otimes (\mathbf{v}-\mathbf{u})\, \G[f] \,d\mathbf{v} \,
=\, \rho \,\mathcal{T}.&
\end{array}\right.
$$
This implies that this collision operator does indeed conserve mass, momentum and energy as imposed.

In this section, we will first recall the  implicit-explicit (IMEX) scheme to the ES-BGK equation proposed in~\cite{bibBGK}. Then we apply our ILW procedure to treat the boundary condition for ES-BGK model case.
\subsection{An IMEX scheme to the ES-BGK equation}
We now introduce the time discretization for the ES-BGK equation~\eqref{eq:boltzmann},~\eqref{eq:ESBGKoperator}
\begin{equation} 
\label{eq:deb}
\left\{
\begin{array}{l}
  \displaystyle{\frac{\partial f}{\partial t}  \,+\, \mathbf{v}\,\cdot\,\nabla_{\mathbf{x}} f \,=\,  
\frac{\tau}{\varepsilon}\,(\mathcal{G}[f] -f),  
   \quad \mathbf{x} \in \Omega\subset\mathbb{R}^{d_{\mathbf{x}}},\, \mathbf{v}\in \mathbb{R}^{3},} 
  \\
  \,
  \\
  f(0,\mathbf{x},\mathbf{v})  \,=\, f_{0}(\mathbf{x},\mathbf{v}), \quad \mathbf{x}\in\Omega,\,\mathbf{v}\in\mathbb{R}^{3}, 
\end{array}\right.
\end{equation}
where $\tau$ depends on $\rho$, $\mathbf{u}$ and $T$. 

The time discretization is an IMEX scheme.
 Since the convection term in (\ref{eq:deb}) is not stiff, we will
treat it explicitly. The source terms on the right hand side of 
(\ref{eq:deb}) will be handled using an implicit solver. We simply apply a first order IMEX scheme, 
\begin{equation} 
\label{sch:perturb}
\left\{
\begin{array}{l}
  \displaystyle{\frac{f^{n+1}-f^n}{\Delta t }  + \mathbf{v}\cdot\nabla_{\mathbf{x}} f^n \,=\, 
\frac{\tau^{n+1}}{\varepsilon}(\mathcal{G}[f^{n+1}]-f^{n+1}),}
  \\
  \,
  \\
  f^0(\mathbf{x},\mathbf{v})  = f_{0}(\mathbf{x},\mathbf{v})\,.
\end{array}\right.
\end{equation}
This can be written as
\begin{eqnarray}
\label{AP-1}
f^{n+1} &=& \frac{\varepsilon}{\varepsilon+\tau^{n+1}\Delta t} \left[f^n - \Delta t  \,\mathbf{v}\cdot\nabla_{\mathbf{x}} f^n\right] \,+\, \frac{\tau^{n+1} \Delta t}{\varepsilon+\tau^{n+1} \Delta t}\,\mathcal{G}[f^{n+1}],
 \end{eqnarray}
where $\mathcal{G}(f^{n+1})$ is the anisotropic Maxwellian distribution computed from $f^{n+1}$. Although (\ref{AP-1}) appears  nonlinearly implicit, since the computation of $f^{n+1}$ requires the knowledge of $\mathcal{G}[f^{n+1}]$, it can be solved explicitly. Specifically, upon multiplying (\ref{AP-1}) by $\phi(\mathbf{v})$ 
defined by
$$
\phi(\mathbf{v}) \,:=\, \left(1,\mathbf{v},\frac{|\mathbf{v}|^2}{2}\right) 
$$
and use the conservation properties of $\mathcal{Q}$ and the definition of $\mathcal{G}[f]$ in~\eqref{eq:macro}, we define the macroscopic quantity $U$ by  $U:=(\rho,\rho\,\mathbf{u},E)$ computed from $f$ and get
$$
U^{n+1} = \frac{\varepsilon}{\varepsilon+\tau^{n+1} \Delta t} \int_{\mathbb{R}^{3}} \phi(\mathbf{v})\,(f^n-\Delta t\,  \mathbf{v} \cdot \nabla_{\mathbf{x}} f^n) \, d\mathbf{v} 
+ \frac{\tau^{n+1}\Delta t}{\varepsilon+\tau^{n+1} \Delta t} \,\int_{\mathbb{R}^{3}} \phi(\mathbf{v})\mathcal{G}[f^{n+1}](\mathbf{v})\,d\mathbf{v},
$$
or simply
\begin{equation}
\label{AP-11}
U^{n+1} =  \int_{\mathbb{R}^{3}}  \phi(\mathbf{v})\,(f^n-\Delta t \mathbf{v} \cdot \nabla_{\mathbf{x}} f^n)  \, d\mathbf{v} \,.
\end{equation}
Thus $U^{n+1}$ can be obtained explicitly. This gives $\rho^{n+1},
\mathbf{u}^{n+1}$ and $T^{n+1}$.
Unfortunately, it is not enough to define $\mathcal{G}[f^{n+1}]$ for which we need $\rho^{n+1}\,\Theta^{n+1}$. Therefore, we define the tensor $\Sigma$ by
\begin{equation}
\label{sigm}
\Sigma^{n+1} := \int_{\mathbb{R}^{3}} \mathbf{v} \otimes \mathbf{v}\, f^{n+1}\,d\mathbf{v} \,= \;\rho^{n+1}\, \left(\Theta^{n+1} \,+\, \mathbf{u}^{n+1} \otimes \mathbf{u}^{n+1}\right)
\end{equation}
and multiply the scheme (\ref{AP-1}) by $\mathbf{v} \otimes \mathbf{v}$. Using the fact that 
$$
\int_{\mathbb{R}^{3}} \mathbf{v} \otimes \mathbf{v}\, \mathcal{G}[f](\mathbf{v})\,d\mathbf{v} =  \rho\, \left(\mathcal{T}\,+\, \mathbf{u} \otimes \mathbf{u}\right),
$$
and (\ref{sigm}), we get that
\begin{eqnarray}
\label{AP-111}
\Sigma^{n+1} &=& \frac{\varepsilon}{\varepsilon+(1-\nu)\,\tau^{n+1}\,\Delta t}\left(\Sigma^n \,-\, \Delta t\,\int_{\mathbb{R}^{3}} \mathbf{v} \otimes \mathbf{v} \,\mathbf{v} \cdot \nabla_{\mathbf{x}} f^n d\mathbf{v} \right)
\\
\,
\nonumber
\\
&& +\, \frac{(1-\nu)\,\tau^{n+1}\,\Delta t}{\varepsilon+(1-\nu)\,\tau^{n+1}\,\Delta t} \, \rho^{n+1}\,\left(T^{n+1} \,{\rm I} \,+\, \mathbf{u}^{n+1}\otimes \mathbf{u}^{n+1} \right)   
\nonumber
\end{eqnarray}
Now $\mathcal{G}[f^{n+1}]$ can be obtained explicitly from $U^{n+1}$ and $\Sigma^{n+1}$ and then $f^{n+1}$ from  (\ref{AP-1}).

Finally the scheme reads
\begin{equation}
\label{AP-final}
\left\{
\begin{array}{lll}
U^{n+1} &=&  \displaystyle{\int_{\mathbb{R}^{3}}  \phi(\mathbf{v})\,(f^n-\Delta t \mathbf{v} \cdot \nabla_{\mathbf{x}} f^n)  \, d\mathbf{v},} 
\\
\,
\\
\Sigma^{n+1} &=& \displaystyle{\frac{\varepsilon}{\varepsilon+(1-\nu)\,\tau^{n+1}\,\Delta t}\left(\Sigma^n \,-\, \Delta t\,\int_{\mathbb{R}^{3}} \mathbf{v} \otimes \mathbf{v} \,\mathbf{v} \cdot \nabla_{\mathbf{x}} f^n d\mathbf{v} \right)}
\\
\,
\\
&&\displaystyle{ +\, \frac{(1-\nu)\,\tau^{n+1}\,\Delta t}{\varepsilon+(1-\nu)\,\tau^{n+1}\,\Delta t} \, \rho^{n+1}\,\left(T^{n+1} \,{\rm I} \,+\, \mathbf{u}^{n+1}\otimes \mathbf{u}^{n+1} \right).}
\\
\,
\\
f^{n+1} &=& \displaystyle{\frac{\varepsilon}{\varepsilon+\tau^{n+1}\Delta t} \left[f^n - \Delta t  \,\mathbf{v}\cdot\nabla_{\mathbf{x}} f^n\right] \,+\, \frac{\tau^{n+1} \Delta t}{\varepsilon+\tau^{n+1} \Delta t}\,\mathcal{G}[f^{n+1}]},   
\end{array}\right.
\end{equation}

The scheme~\eqref{AP-final} is an AP scheme for~\eqref{sch:perturb}. On the one hand, although~\eqref{sch:perturb} is nonlinearly implicit, is can be solved explicitly. On the other hand, the scheme~\eqref{AP-final} preserves the correct asymptotic~\cite{bibFS2011}, which means when holding the mesh size  and time step fixed and letting the Knudsen number go to zero, the scheme becomes a suitable scheme for the limiting hydrodynamic models.



\subsection{Inverse Lax-Wendroff procedure for boundary conditions}
We have described the numerical method for boundary condition to general kinetic equations in spatially 1D and 2D case. To implement this method, it remains to replace the collision operator $\mathcal{Q}(f)$ in~\eqref{eq:1Dreformulation} or~\eqref{eq:2Dreformulation} by the ES-BGK operator~\eqref{eq:ESBGKoperator}.

Assume that the approximation to the distribution function  at the boundary $f(\mathbf{x}_p,\mathbf{v})$ is known for all $\mathbf{v}\in \mathbf{V}_{h}$. Then, the macroscopic quantities $\rho$, $\mathbf{u}$ and ${T}$ at the boundary point $\mathbf{x}_p$ can be obtained using~\eqref{eq:macro} and ~\eqref{eq:ESBGKmacro}. Therefore,  substituting these macroscopic quantities in~\eqref{eq:stress}, we compute the stress tensor $\Theta$ at the boundary point $\mathbf{x}_p$, such that the corrected tensor $\mathcal{T}(\mathbf{x}_p)$. Thus $\mathcal{G}[f]$ is computed for all points $(\mathbf{x}_p,\mathbf{v})$, where $\mathbf{v}\in \mathbf{V}_{h}$.


\section{Numerical examples}
\label{sec:Num}
\setcounter{equation}{0}
In this section, we present a large variety of test cases in $1d_x$ and $2d_x$  in space and three dimensional in velocity space showing the effectiveness of our method to get an accurate solution of Boltzmann type equations set in a complex geometry with different boundary conditions.  We first give an example on a flow generated by gradients of temperature, which has already been treated by DSMC or other various methods \cite{bibF}. 

Finally, we present some numerical results in $2d_x$.


\subsection{Smooth solutions}
We consider the ES-BGK equation (\ref{eq:boltzmann})-(\ref{eq:ESBGKoperator})
$$
\left\{
\begin{array}{l}
\displaystyle \frac{\partial f}{\partial t}  + v_x \frac{\partial f}{\partial x} \,=\, \mathcal{Q}(f), , \quad x\in (-0.5,\,0.5), \, v\in\RR^3,
\\
\,
\\
f(t=0) = f_0,
\end{array}\right.
$$
with an initial datum $f_0$ which is a perturbation of the constant state in space and a Maxwellian distribution function in velocity, that is,
$$
f_0(x,v)= \frac{\rho_0(x)}{(2\pi)} \,\exp\left(-\frac{|v|^2}{2}\right), \quad x\in (-0.5,\,0.5), \, v\in\RR^3
$$
with a density $\rho_0= 1+\,0.1\,\cos(2\pi\,x)$. We consider purely diffusive boundary conditions with a wall temperature $T_w=1$.  The solution is expected to be smooth for short time and then may develop a discontinuity at the boundary, which may propagate in the physical domain. 

We perform several numerical simulations on a time interval $[0,t_{end}]$ with $t_{end}=1$,  a computational domain in space $I_0=[-\pi/6,\pi/6]$ such that $(-1/2,1/2)\subset I_0$ and a domain in velocity  $\mathbf{V}=[-8,8]^3$. Then, we choose a grid in space for $I_0$ constituted of $n_x=n$ points and a grid $\mathbf{V}_{h}$ for the velocity space with $n_v=n$ points for each direction with respectively $n=32$, $64$,..., $n=512$. Let us emphasize that the boundary points $x=-1/2$ and $x=1/2$ are not exactly located on a grid point. Since we don't know an exact solution of the problem, we compute relative errors. More precisely, an estimation of the relative error in $L^{1}$ norm at time $T$ is given by
\begin{equation*}
e_{2 h}=\Vert f_{h}(T)-f_{2h}(T)\Vert_{L^{1}},
\end{equation*}
where $f_{h}$ represents the approximation computed from a mesh of size $h=(\Delta x,\Delta v)$. The numerical scheme is said to be $k$-th order if $e_{2h} \leq C \|h\|^{k}$, for all $0<\|h\| \ll 1$.

In Table \ref{table1} we compute the order of convergence in  $L^{1}$ norm  of our numerical methods. We can clearly see the expected second order convergence. Moreover, we verify experimentally that our scheme is also second-order accurate at the boundary since the discontinuity occuring at $v_x=0$ is perfectly located.

\begin{table}[!ht]
\centering
\begin{tabular}{|c|c|c|c|c|}
\hline $n_{x}\times n_v$  & $L^{1}$ error  & Order & $L^{1}$ error at the boundary & Order \\ 
\hline 
       $32^2$     & $ 8.8833 \,10^{-4} $  & X &  $ 3.909 \,10^{-3} $  & X \\ 
\hline
       $64^2$     & $ 2.5221 \,10^{-4} $  & 1.94 &  $ 5.832 \,10^{-3} $  & X  \\ 
\hline
       $128^2$    & $ 6.5511 \,10^{-5} $  & 1.88 &  $ 2.341 \;10^{-4} $  & 4.1 \\
\hline
       $256^2$    & $ 1.7829 \, 10^{-5} $  &1.91 &  $ 5.811 \,10^{-5} $  & 2.01 \\
\hline
       $512^2$    & $ 4.4571 \, 10^{-6} $  & 2 &    $ 1.573 \,10^{-5} $  & 1.89 \\
    \hline
\end{tabular} 
\caption{Smooth solutions:  {\em  Experimental order of convergence in $L^{1}$ norm.}}
\label{table1}
\end{table}

\begin{figure}[htbp]
\begin{tabular}{cc}
\includegraphics[width=7.5cm]{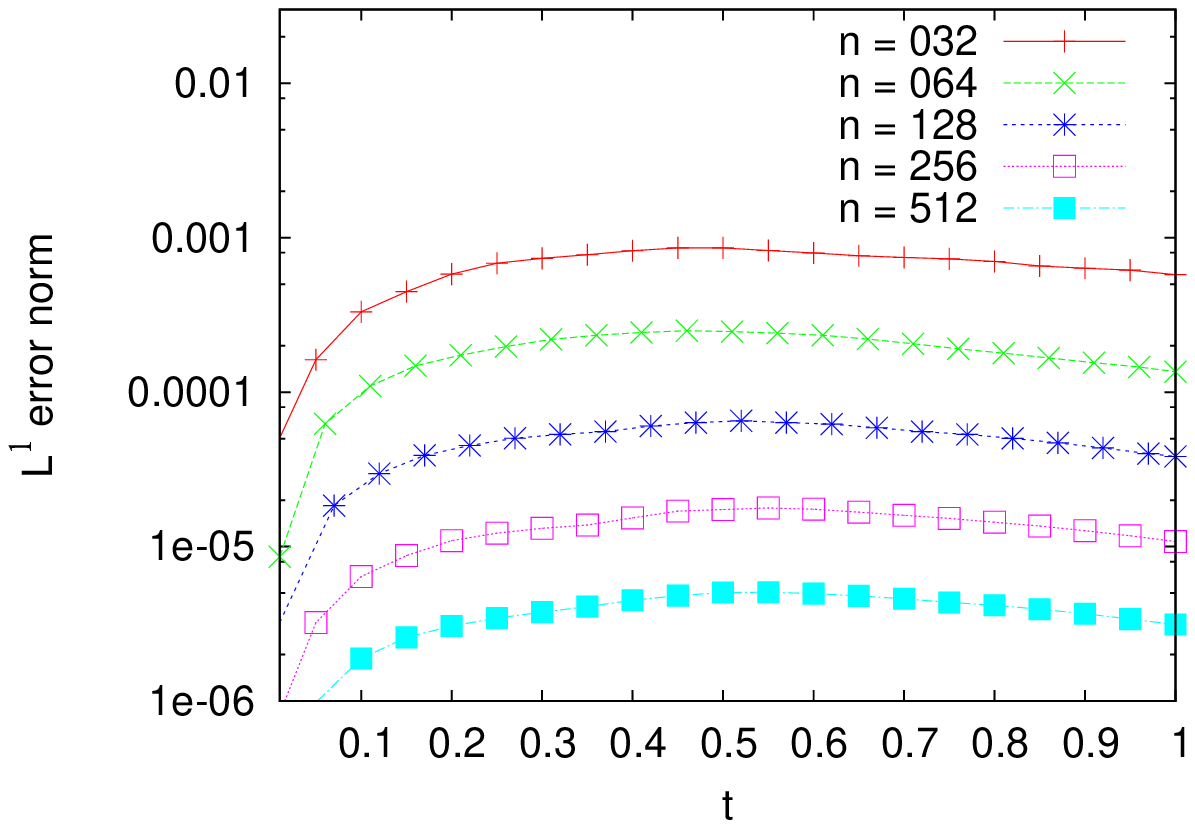} &
\includegraphics[width=7.5cm]{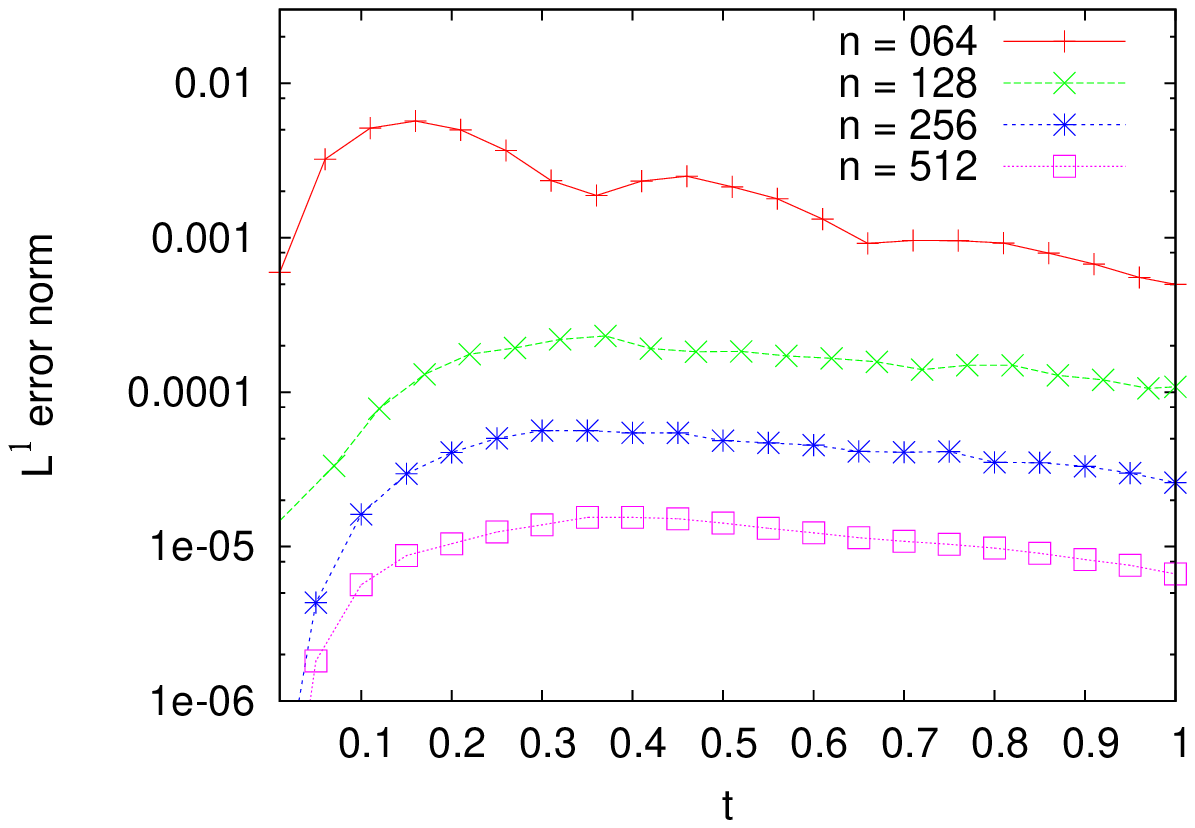} 
\\
(1)&(2)
\end{tabular}
\caption{Smooth solutions: {\em Experimental order of convergence in $L^{1}$ norm (1) in the physical domain (2) at the boundary.} }
\label{Fig_accur}
\end{figure}

\subsection{Flow generated by a gradient of temperature}
\label{sec4-2}
We consider the ES-BGK equation (\ref{eq:boltzmann})-(\ref{eq:ESBGKoperator}), 
$$
\left\{
\begin{array}{l}
\displaystyle \frac{\partial f}{\partial t}  + v_x \frac{\partial f}{\partial x} \,=\, \frac{1}{\varepsilon}\,\mathcal{Q}(f),\,\, x\in (-1/2,1/2),\,v\in\RR^3,
\\
\,
\\
\displaystyle f(t=0,x,v) = \frac{1}{2\pi \, T_0} \,\exp\left(-\frac{|v|^2}{2\,T_0}\right),  
\end{array}\right.
$$
with  $T_0(x)= 1$ and we assume purely diffusive boundary conditions on $x=-1/2$ and $x=1/2$, which can be written as 
$$
f(t,x,v) = \mu(t,x)\; f_w(v), \,{\rm if }\,(x,v_x)\in\{-1/2\}\times\RR^+\,{\rm and }\,(x,v_x)\in\{1/2\}\times\RR^-,
$$
where $\mu$ is given by (\ref{eq:MU}). This problem has already been studied in \cite{wagner} using DSMC for the Boltzmann equation or using deterministic approximation using a BGK model for the Boltzmann equation in \cite{Aoki94,bibF}. 

Here we apply our numerical scheme with the ES-BGK operator (\ref{eq:ESBGKoperator}) and choose a computational domain in space $I_0=[-\pi/6,\pi/6]$ such that $(-1/2,1/2)\subset I_0$ and $[-8,8]^3$ for the velocity space with a number grid points  $n=32$ in each direction and the time step $\Delta t=0.001$.

The main issue here is to capture the correct steady state for which the pressure is a perturbation of a constant state with a Knudsen layer at the boundary \cite{Aoki94,wagner}.  

\begin{figure}[htbp]
\begin{tabular}{cc}
\includegraphics[width=7.5cm]{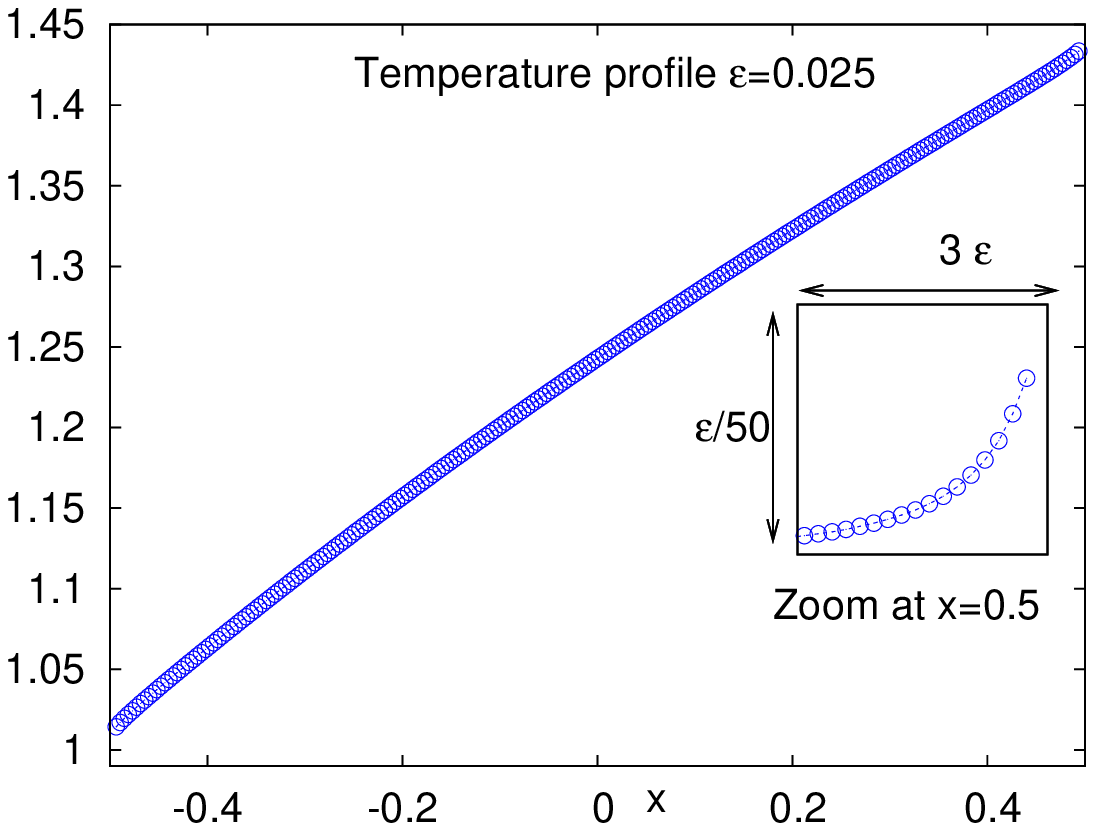} &
\includegraphics[width=7.5cm]{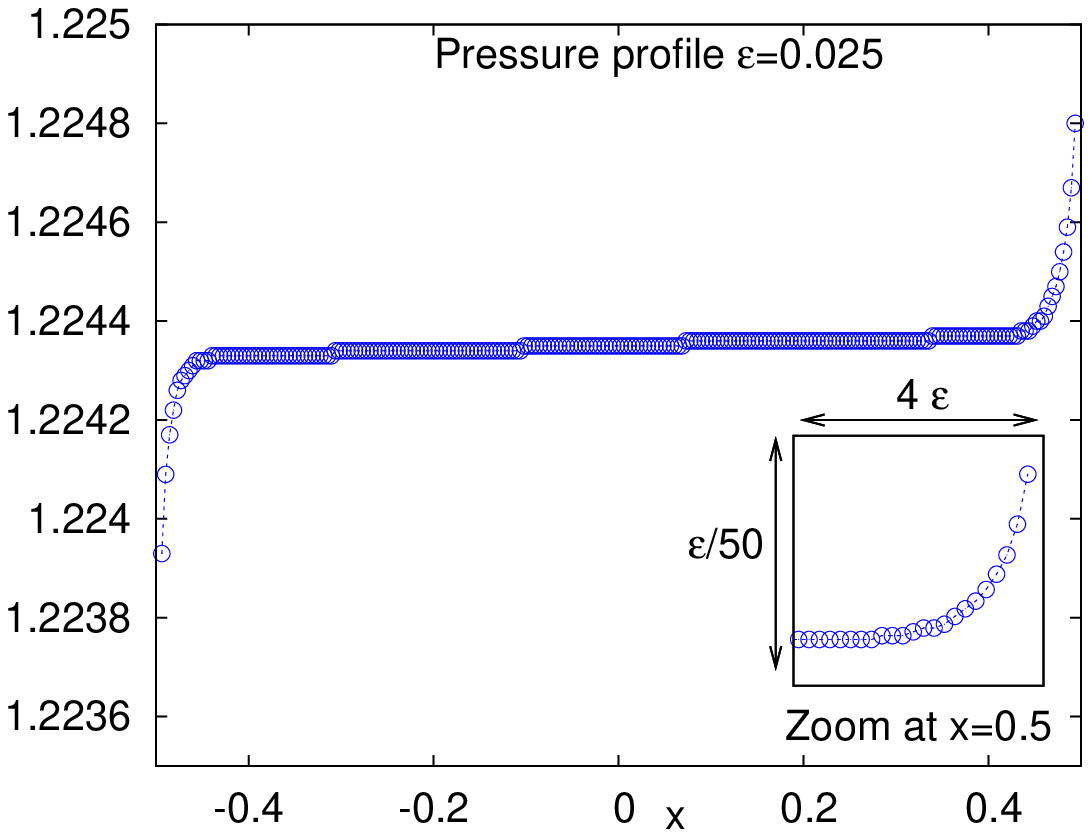} 
\\
\includegraphics[width=7.5cm]{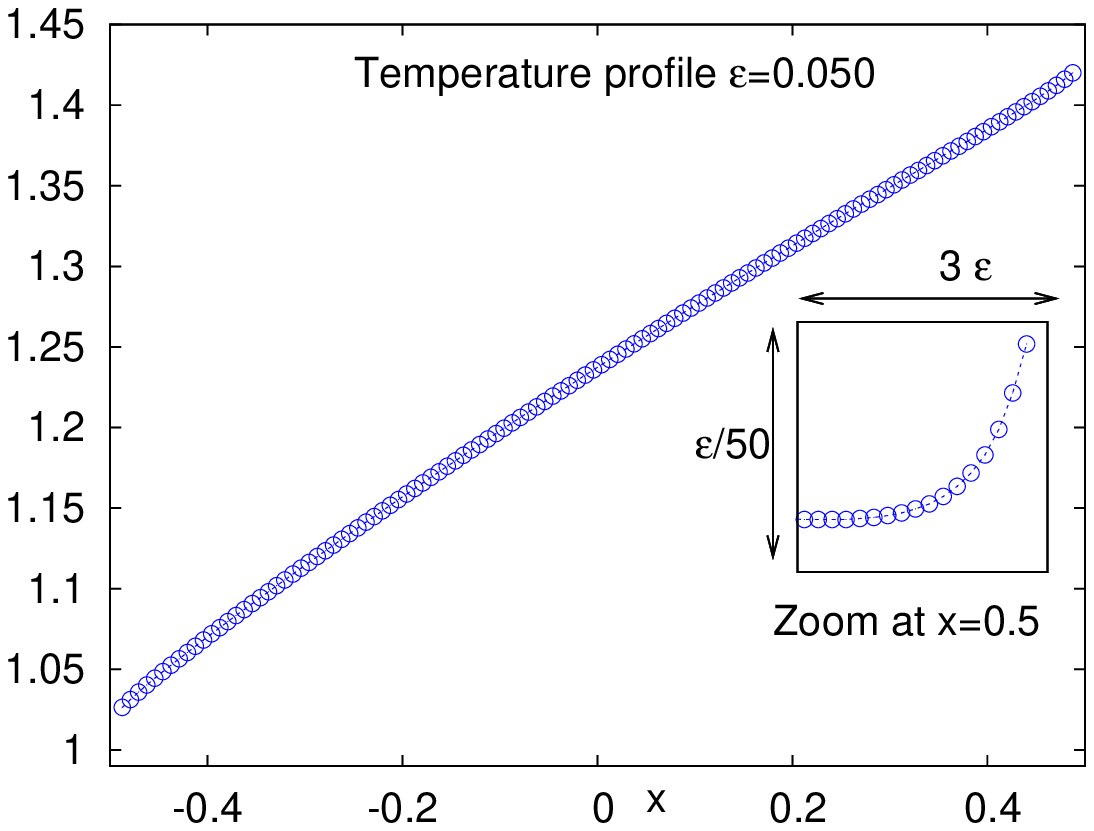} &
\includegraphics[width=7.5cm]{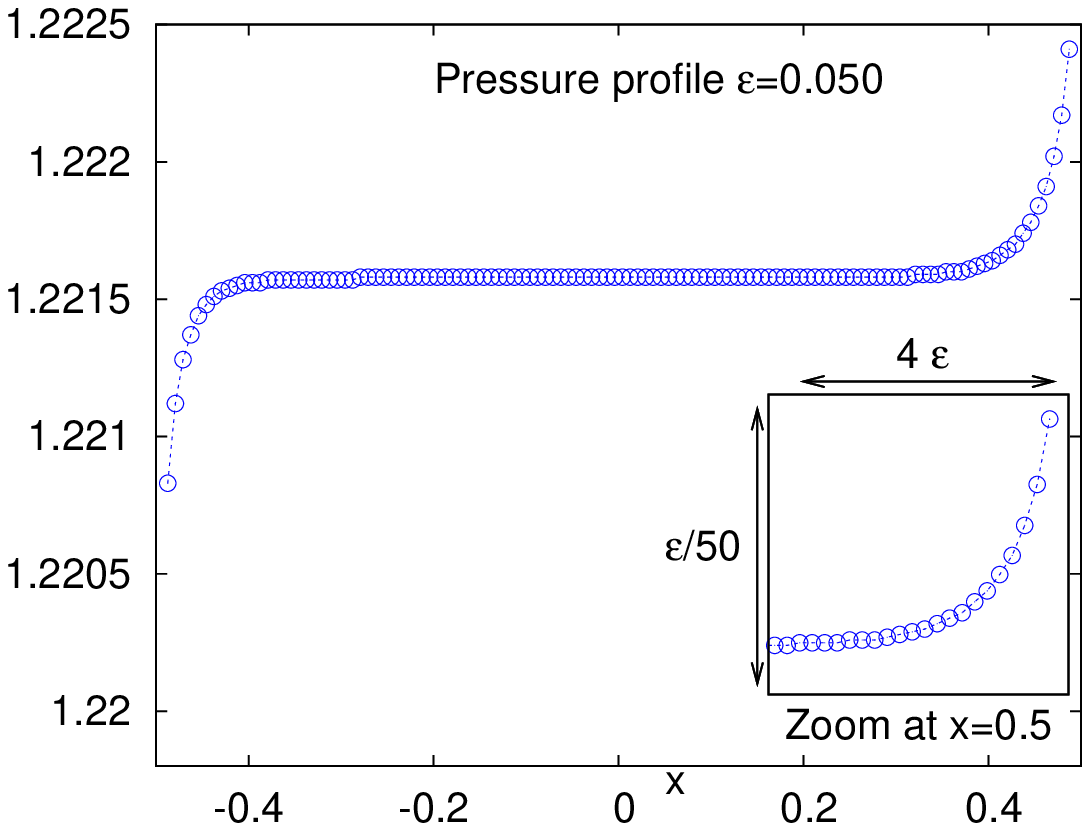}
\\
\includegraphics[width=7.50cm]{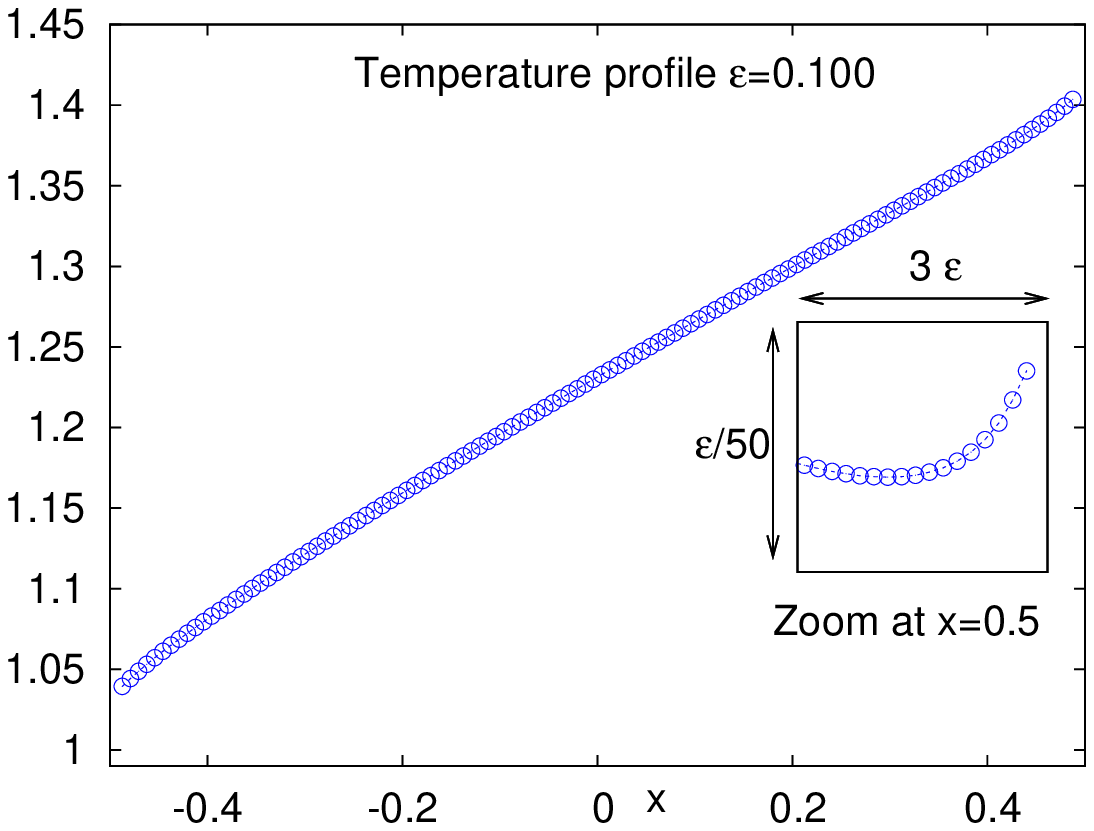} &
\includegraphics[width=7.50cm]{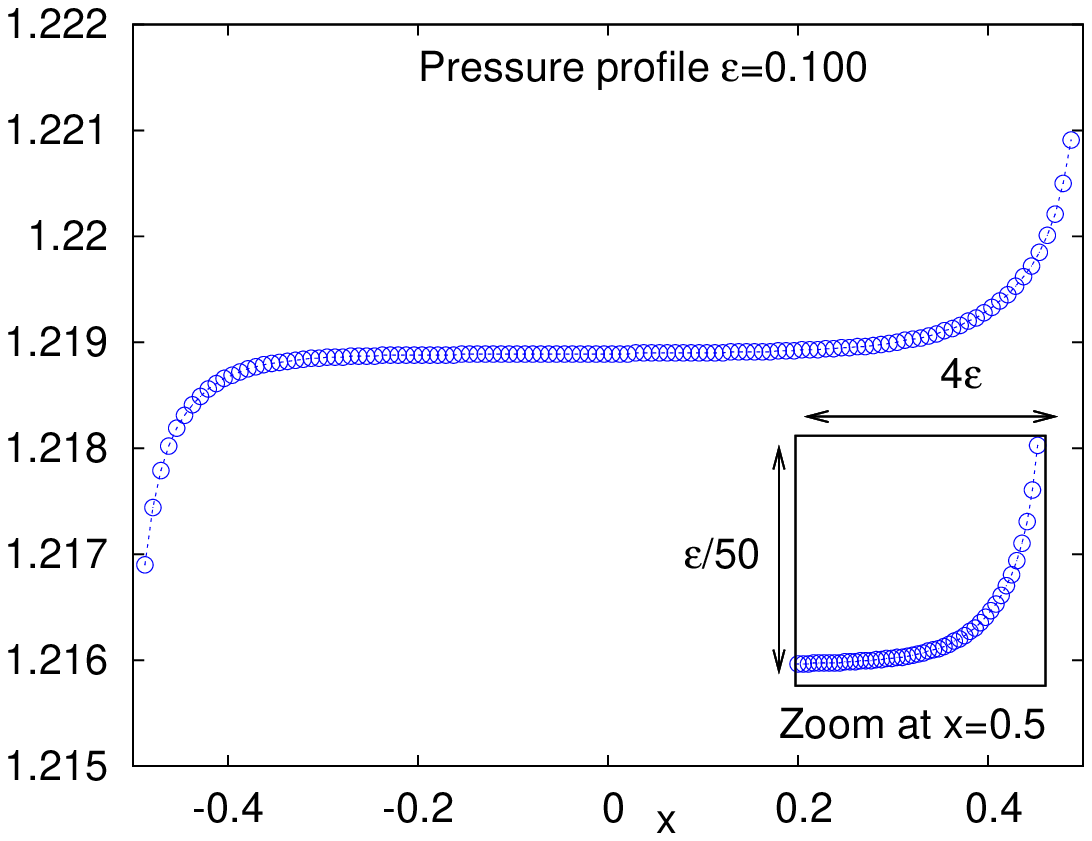}     
\\
(1)&(2)
\end{tabular}
\caption{Flow generated by a gradient of temperature: {\em (1) temperature (2)  pressure   for various Knudsen numbers $\varepsilon = 0.025$, $0.05$ and $0.1$.}}
\label{Fig_bl_1}
\end{figure}

In Figure~\ref{Fig_bl_1}, we represent the stationary solution (obtained approximately at time $t_{end}=25$ for $\varepsilon=0.1$ up to $t_{end}=75$ for $\varepsilon=0.025$) of the temperature and the pressure profile. The results are in a qualitative good agreement with those already obtained in \cite{wagner} with DSMC. More precisely, the boundary layer (Knudsen layer) appears in the density and temperature as well as the pressure, but it is small for all the quantities. The magnitude in the dimensionless density, temperature, and pressure is of order of $\varepsilon$ and the thickness of the layer is, say $O(\varepsilon)$. In the density and temperature profiles, we cannot observe it unless we magnify the profile in the vicinity of the boundary (see the zoom in Figure~\ref{Fig_bl_1}). Instead, since the pressure is almost constant in the bulk of the gas, we can observe perfectly the boundary layer by magnifying the entire profile. Let us emphasize that, as it is shown in Figure~\ref{Fig_bl_1} the Knudsen layer is a kinetic effect, which disappears in the fluid limit ($\varepsilon \rightarrow 0$). 

These results provide strong evidence that the present treatment of boundary conditions using WENO extrapolation and inverse Lax-Wendroff  method can be used to determine the state of a gas under highly non-equilibrium conditions. Using deterministic methods, we can investigate the behavior of gases for situations in which molecular diffusion is important e.g., thermal diffusion. 

Also let us mention that a quantitative comparison between our results ($3d_{\mathbf{v}}$ with ES-BGK operator) and  \cite{wagner} ($3d_{\mathbf{v}}$ Boltzmann with hard sphere potential) or \cite{Aoki94} ($3d_{\mathbf{v}}$ BGK) gives a very good agreement on the values of the Knudsen layer and the values of the pressure inside the domain.

\subsection{High-speed flow through a trapezoidal channel}
\label{sec4-3}
\begin{figure}
\begin{center}
  \begin{tikzpicture}
    \draw   (0,0) node{\includegraphics[width = 12cm]{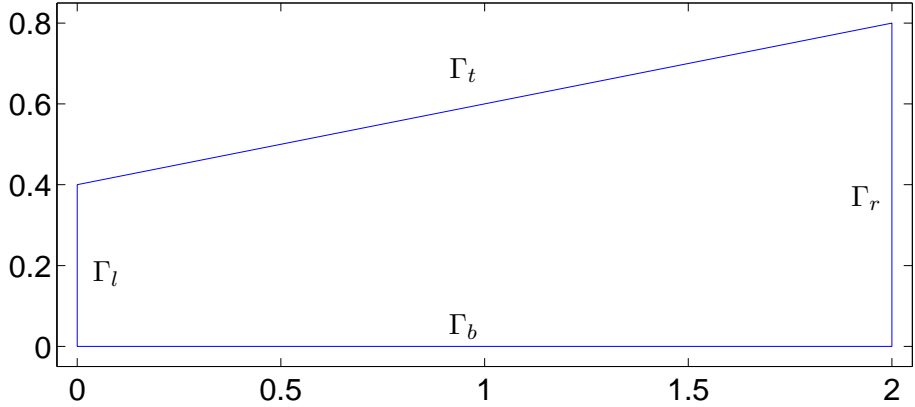}};
    \draw   (-4.7,-1) node{$\Gamma_l$};
    \draw   (5.3,0) node{$\Gamma_r$};
    \draw   (0,1.7) node{$\Gamma_t$};
    \draw   (0,-1.7) node{$\Gamma_b$};
  \end{tikzpicture}
\end{center}
\caption{\label{fig:trapeze}Trapezoidal domain $\Omega_{\mathbf{x}}$.}
\end{figure}

In this section we deal with spatially two-dimensional ES-BGK model in a trapezoidal domain. We attempt to get some steady state as
\begin{equation*}
   \displaystyle  v_x\frac{\partial f}{\partial x}+v_y\frac{\partial f}{\partial y}=\frac{1}{\varepsilon}\mathcal{Q}(f),
\end{equation*}
where $\mathbf{v}\in \Omega_{\mathbf{x}}$ and $\mathbf{v}\in \mathbb{R}^3$. Here we will reproduce  a numerical test performed in~\cite{bibRW} but with our ILW method. The computational domain is a trapezoid
\begin{equation*}
  \Omega_{\mathbf{x}}=\{\mathbf{x}=(x,y),\,\,0<x<a,\,\,0<y<b+x\tan(\alpha)\}
\end{equation*}
as shown in Figure~\ref{fig:trapeze} for the parameters
\begin{equation*}
  a=2.0,\,\,b=0.4,\,\,\alpha=\arctan(0.2).
\end{equation*}
 Boundary conditions are defined separately for each of the four straight pieces
\begin{equation*}
  \partial \Omega_{\mathbf{x}}=\Gamma_l\cup\Gamma_b\cup\Gamma_r\cup\Gamma_t
\end{equation*}
denoting the left, bottom, right and top parts of the boundary respectively.
The bottom part represents the axis of symmetry, so we use specular reflection~\eqref{op:BC} there, {\it i.e.}
\begin{equation*}
  f(\mathbf{x},\mathbf{v})=(f(\mathbf{x},\mathbf{v}-2(\mathbf{v}\cdot\mathbf{n}(\mathbf{x}))\mathbf{n}(\mathbf{x})),\,\,\,\mathbf{x}\in\Gamma_b,\,\,v_y>0.
\end{equation*}
On the right part we are modeling outflow (particles are permanently absorbed), {\it i.e.}
\begin{equation}
  f(\mathbf{x},\mathbf{v})=0,\,\,\,\mathbf{x}\in\Gamma_r,\,\,v_x<0.
  \label{eq:inflow_right}
\end{equation}
On the left part there is an incoming flux of particles, {\it i.e.}
\begin{equation}
  f(\mathbf{x},\mathbf{v})=f_{\text{in}}(\mathbf{x},\mathbf{v})=M_{\text{in}}(\mathbf{v}),\,\,\,\mathbf{x}\in\Gamma_l,\,\,v_x>0,
  \label{eq:inflow_left}
\end{equation}
with an inflow Maxwellian
\begin{equation*}
  M_{\text{in}}(\mathbf{v})=\frac{\rho_{\text{in}}}{(2\pi T_{\text{in}})^{3/2}}\exp\left(-\frac{|\mathbf{v}-V_{\text{in}}|^2}{2T_{\text{in}}}\right).
\end{equation*}
On the top part of the boundary, we consider a diffuse reflection~\eqref{op:BC} of particles, with a  Maxwellian distribution function
\begin{equation*}
  M_{\Gamma_t}(\mathbf{v})=\exp\left(-\frac{|\mathbf{v}|^2}{2T_t}\right).
\end{equation*}

In the numerical experiments we assume
\begin{equation*}
\rho_{\text{in}}=1,\,\,\,T_{\text{in}}=1,\,\,\,T_{t}=1.05
\end{equation*}
and consider the inflow velocity in the form
\begin{equation*}
\label{eq:Vin}
  V_{\text{in}}=\text{Mach}_{\text{in}}\sqrt{\gamma T_{\text{in}}}\begin{pmatrix}1\\0\\0\end{pmatrix},
\end{equation*}
where $\text{Mach}_{\text{in}}=5$ and $\gamma=1.4$.

To start the calculation we take an uniform initial solution equal to the values defined by the left boundary condition:
\begin{equation*}
\label{eq:IC:2D}
f_0(\mathbf{x},\mathbf{v})=\frac{\rho_{\text{in}}}{(2\pi T_{\text{in}})^{3/2}}\exp\left(-\frac{|\mathbf{v}-V_{\text{in}}|^2}{2T_{\text{in}}}\right),\,\,\,\mathbf{x}\in \Omega_{\mathbf{x}},\,\,\mathbf{v}\in\mathbb{R}^3.
\end{equation*}

We define the Mach number from the macroscopic quantities, computing the moments of the distribution function with respect to $\mathbf{v}\in\mathbb{R}^3$, by
\begin{equation*}
  \text{Mach}=\frac{|\mathbf{u}|}{\sqrt{\gamma T}},
\end{equation*}
where $c:=\sqrt{\gamma T}$ is the sound speed.

We apply our inverse Lax-Wendroff method to the boundary conditions at $\partial \Omega_{\mathbf{x}}$. More precisely, we extrapolate first the outflow at ghost points corresponding to the four  straight pieces.  Then we impose directly the inflow at the boundaries $\Gamma_l$ and $\Gamma_r$  by~\eqref{eq:inflow_right},~\eqref{eq:inflow_left}, since they are independent of outflow. While the inflow of $\Gamma_b$ and $\Gamma_t$ is computed by specular and diffuse reflection. Finally we use  inverse Lax-Wendroff procedure to compute inflow  at ghost points.

\begin{figure}
  \begin{tabular}{cc} 
    \includegraphics[width = 8cm]{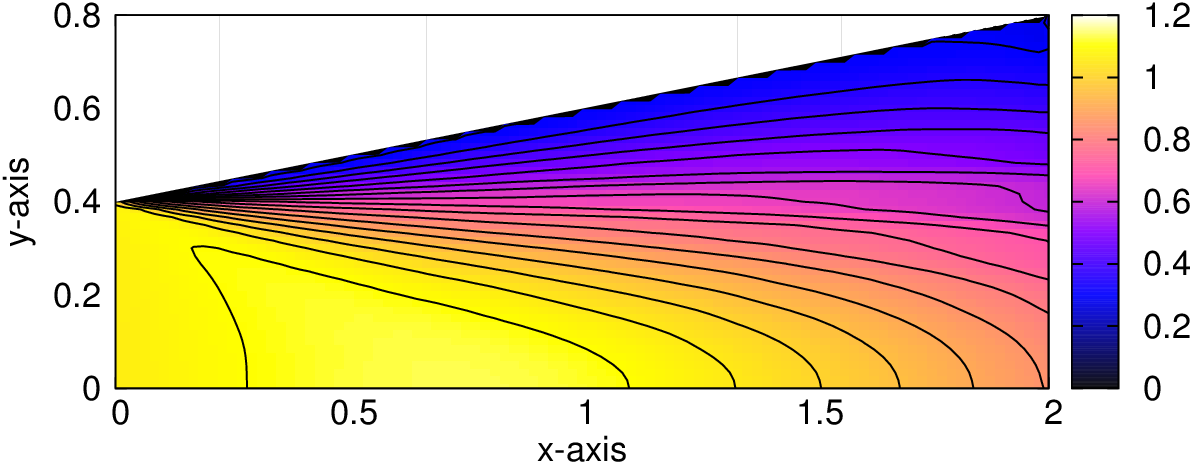}    &
\includegraphics[width = 4.5cm]{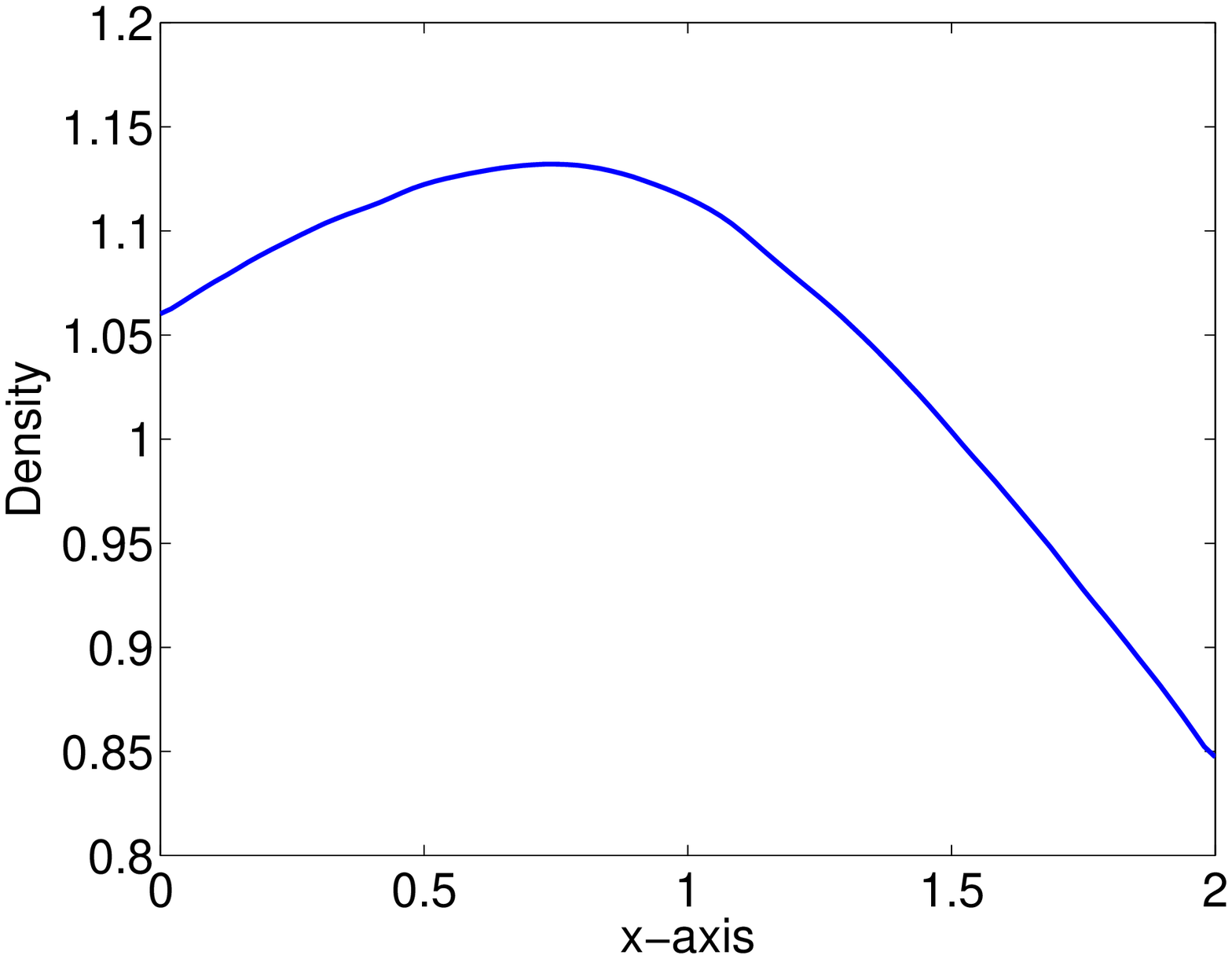}
  \end{tabular}
  \caption{\label{fig:trapeze_density_Kn5}Density with $\varepsilon=5$}
\end{figure}

\begin{figure}
  \begin{tabular}{cc} 
\includegraphics[width = 8cm]{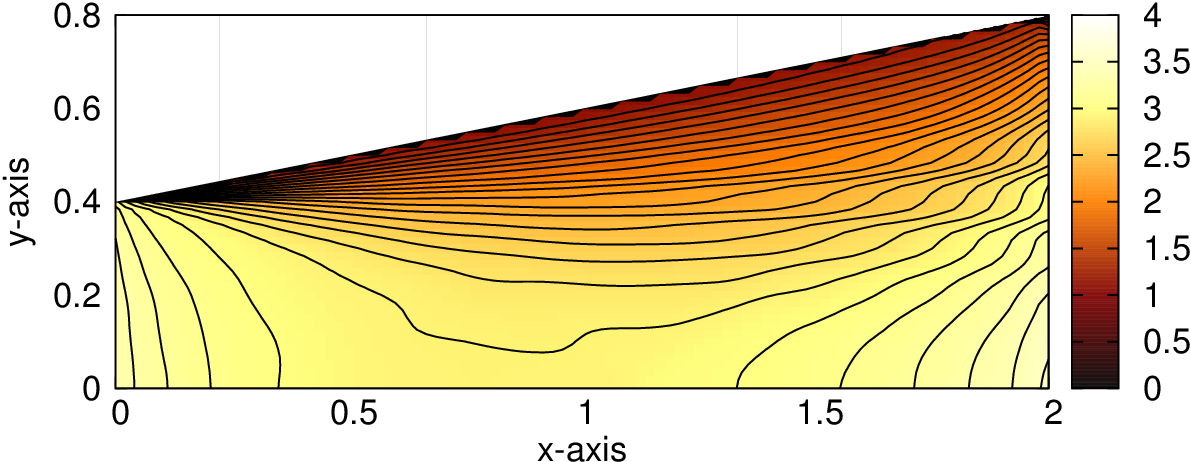}    &
\includegraphics[width = 4.5cm]{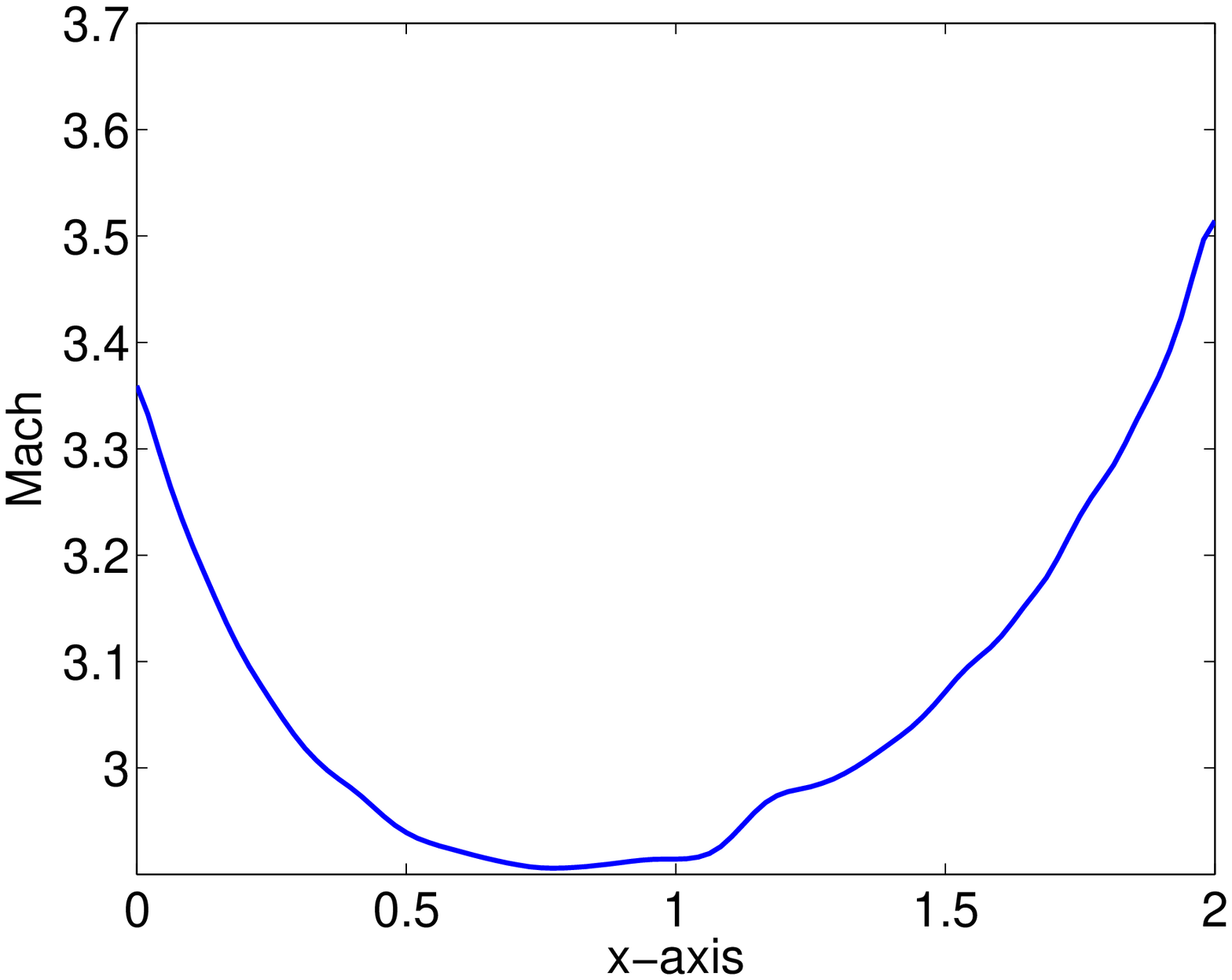}
  \end{tabular}
  \caption{\label{fig:trapeze_mach_Kn5}Mach with $\varepsilon=5$}
\end{figure}

\begin{figure}
  \begin{tabular}{cc} 
    \includegraphics[width = 8cm]{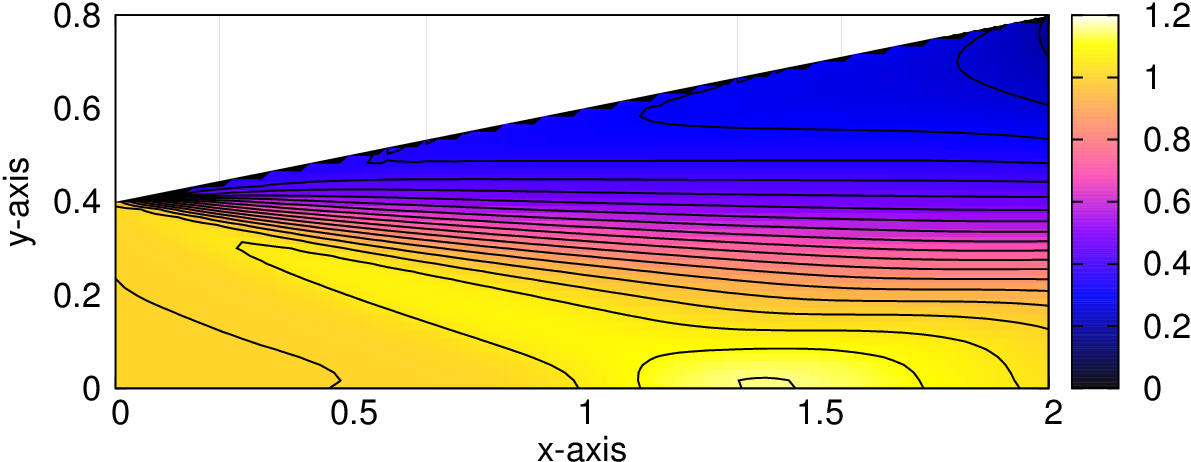}    &
\includegraphics[width = 4.5cm]{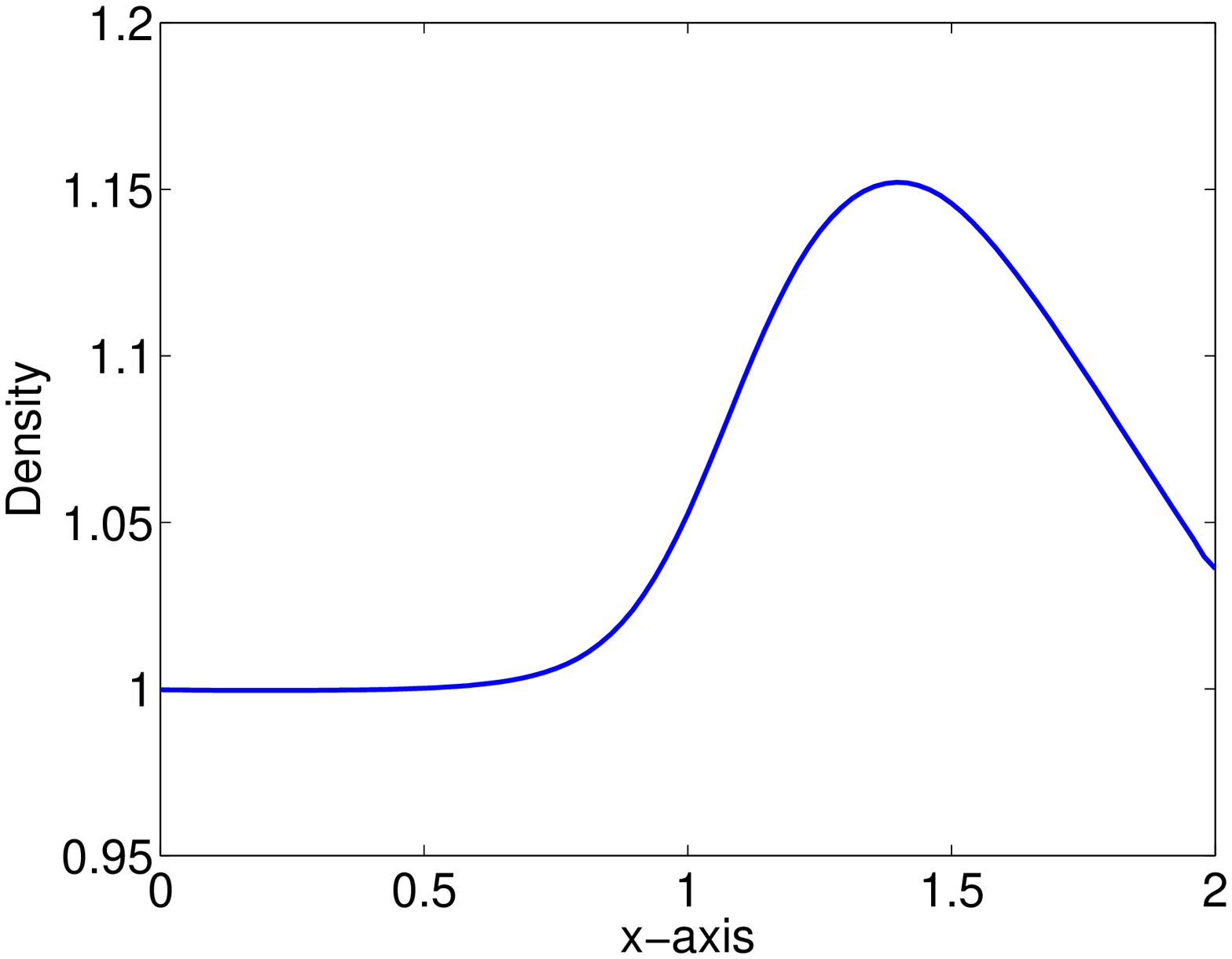}
  \end{tabular}
  \caption{\label{fig:trapeze_density_Kn005}Density with $\varepsilon=0.05$}
\end{figure}

\begin{figure}
  \begin{tabular}{cc} 
    \includegraphics[width = 8cm]{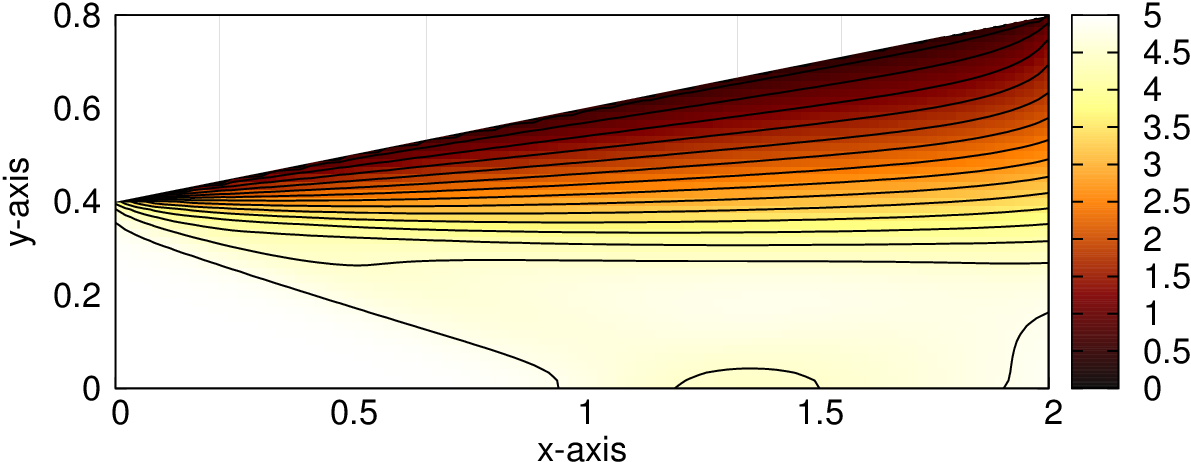}    &
\includegraphics[width = 4.5cm]{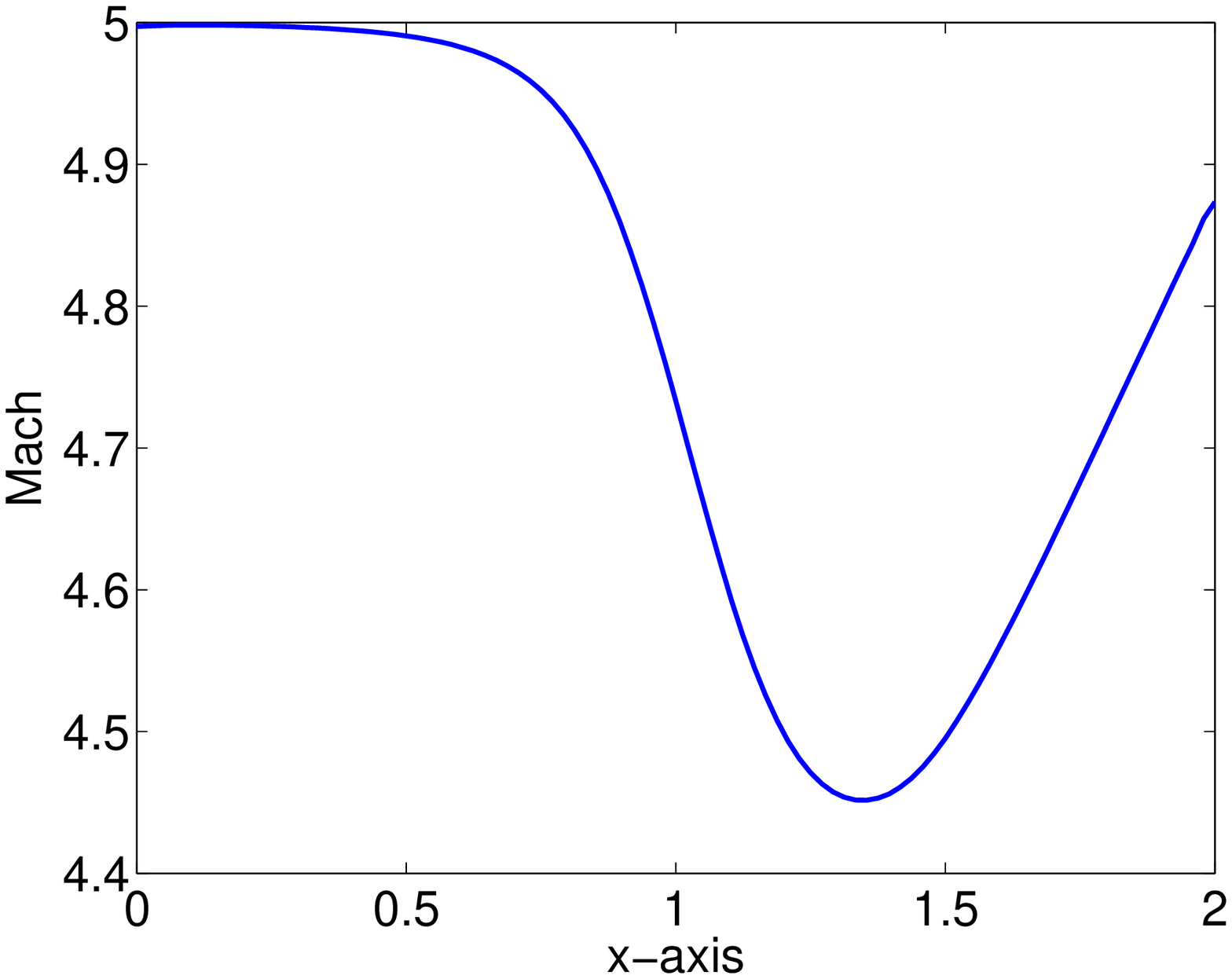}
  \end{tabular}
  \caption{\label{fig:trapeze_mach_Kn005}Mach with $\varepsilon=0.05$}
\end{figure}

In following the sequel, numerical experiments are performed on  a mesh of size $96\times48$   on space domain $\Omega_{\mathbf{x}}$. For velocity space we choose limit domain $[-12,12]\times[-8,8]\times[-8,8]$ with the  grid point number as $64\times48\times12$. Moreover for the ES-BGK operator~\eqref{eq:ESBGKoperator} we choose $\nu=-0.5$. We first consider the weak collision case, {\it i.e.} $\text{Kn}=5$.  In Figures~\ref{fig:trapeze_density_Kn5}--\ref{fig:trapeze_mach_Kn5}, we show on the left hand side, the contour plots of the density, the temperature and the Mach number while the right hand side plots show the absolute values of these quantities plotted along the axis of symmetry $y=0$. We observe that the flow changes when we consider different Knudsen numbers $\text{Kn}=0.05$ and $\text{Kn}=5$. The corresponding results are shown in Figures~\ref{fig:trapeze_density_Kn005}--\ref{fig:trapeze_mach_Kn005}. The significant difference between these two case can be observed in Mach number. In the case $\text{Kn}=5$, the Mach number reaches its maximum at $\Gamma_r$ while in the case $\text{Kn}=0.05$  its maximum is at  $\Gamma_l$. We can observe also in the case $\text{Kn}=0.05$ that there is a clear maximum of the density  in the middle of the domain. In the same region the temperature reaches its maximum.

\begin{figure}
  \begin{tabular}{cc} 
\includegraphics[width = 8cm]{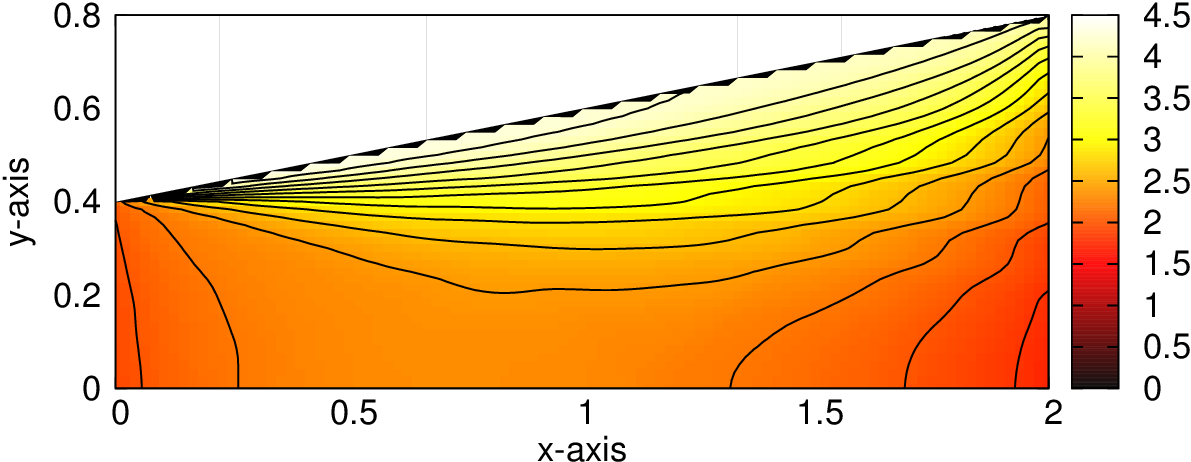}    &
\includegraphics[width = 4.5cm]{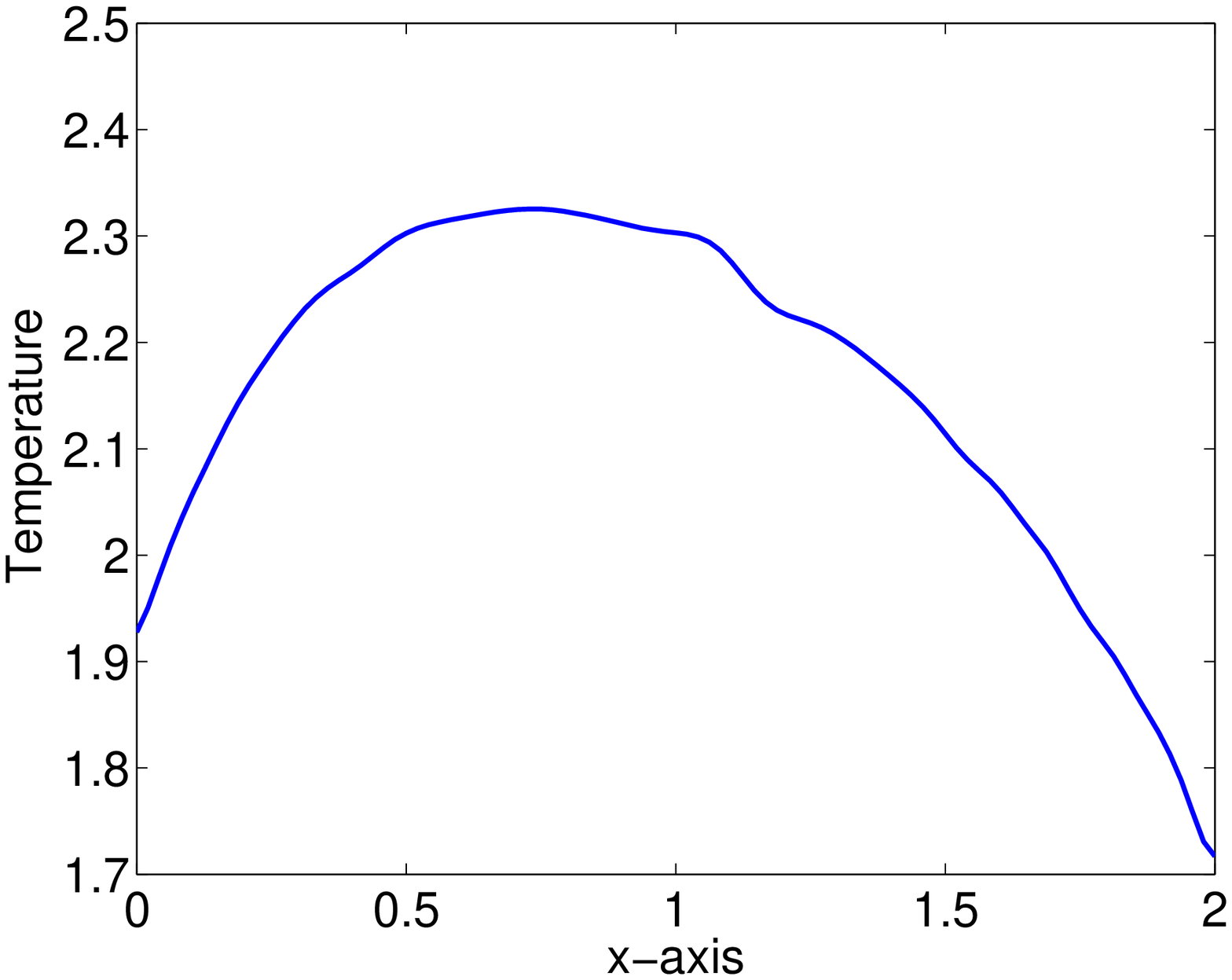}
  \end{tabular}
  \caption{Temperature, $\varepsilon=5$}
\end{figure}

\begin{figure}
  \begin{tabular}{cc} 
    \includegraphics[width = 8cm]{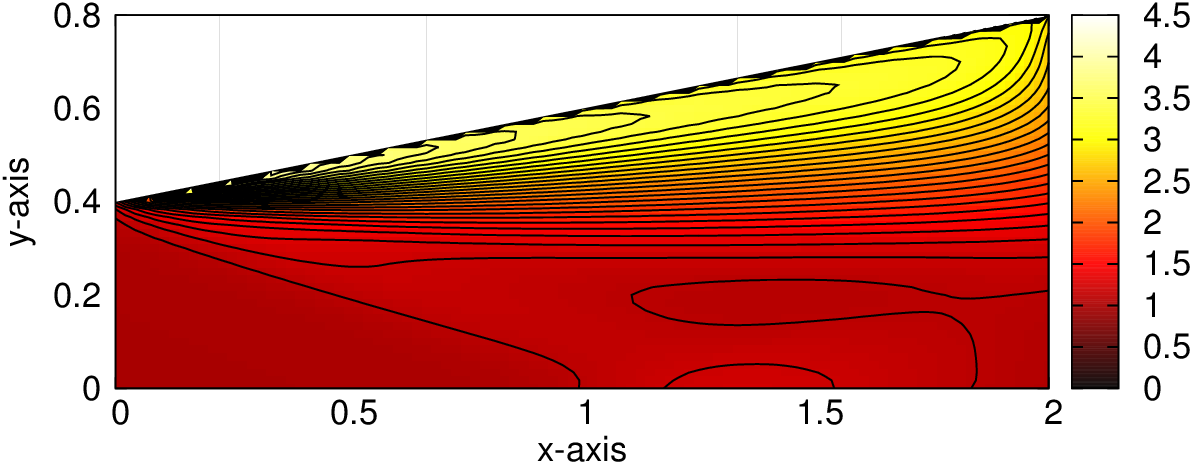}    &
\includegraphics[width = 4.5cm]{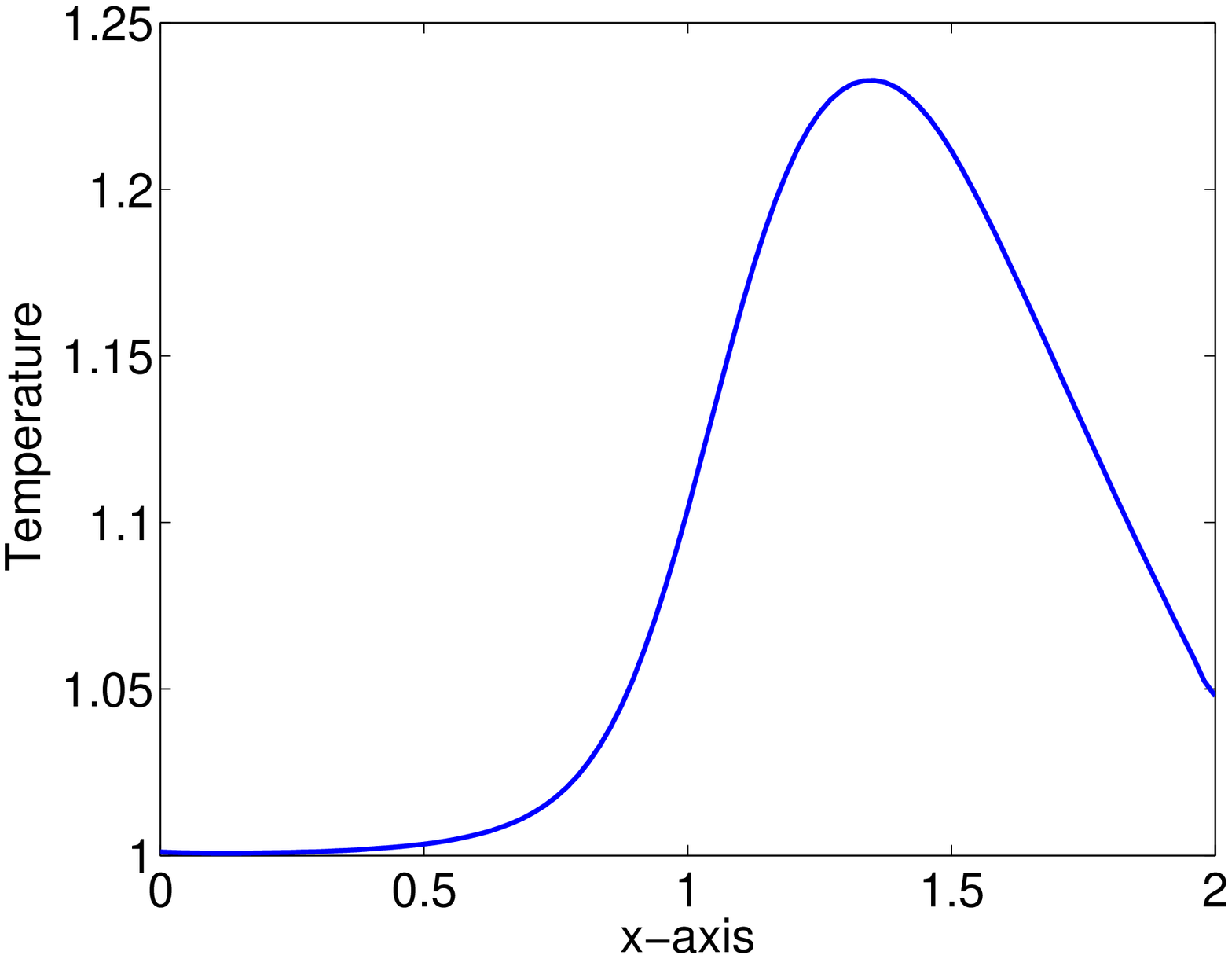}
  \end{tabular}
  \caption{Temperature with $\varepsilon=0.05$}
\end{figure}

\subsection{High-speed flow around an object}
\label{sec4-4}
\begin{figure}
\begin{center}
  \begin{tikzpicture}
    \draw   (0,0) node{\includegraphics[width = 8cm]{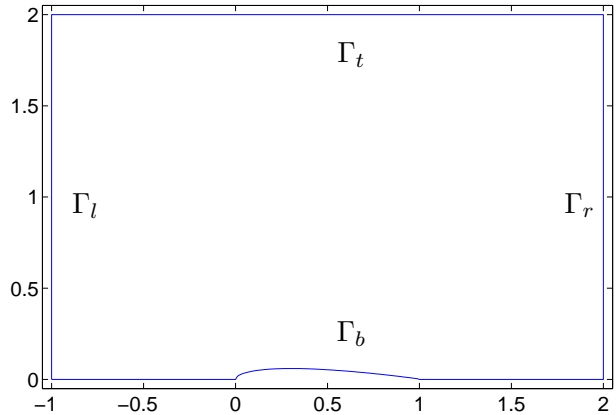}};
    \draw   (-3,0) node{$\Gamma_l$};
    \draw   (3.5,0) node{$\Gamma_r$};
    \draw   (0.5,2.0) node{$\Gamma_t$};  
    \draw   (0.5,-1.7) node{$\Gamma_{b}$};   
  \end{tikzpicture}
\end{center}
\caption{\label{fig:wing}Domain including a half airfoil $\Omega_{\mathbf{x}}$.}
\end{figure}

In this section, we desire to simulate high-speed airflow around a half airfoil  (see Figure~\ref{fig:wing}). The   boundary is separated by four parts
\begin{equation*}
  \partial \Omega_{\mathbf{x}}=\Gamma_l\cup\Gamma_b\cup\Gamma_r\cup\Gamma_t.
\end{equation*}
 On the right and left hand sides $\Gamma_r$, $\Gamma_l$, we use the same incoming flux~\eqref{eq:inflow_right},~\eqref{eq:inflow_left}. On the top  part  $\Gamma_t$, the incoming flow is given by the initial value
\begin{equation*}
  f(\mathbf{x},\mathbf{v})=\frac{\rho_{\text{in}}}{(2\pi T_{\text{in}})^{3/2}}\exp\left(-\frac{|\mathbf{v}-V_{\text{in}}|^2}{2T_{\text{in}}}\right),\,\,\,\mathbf{x}\in \Gamma_t,\,\,\mathbf{v}\in\mathbb{R}^3.
\end{equation*}
Finally at the bottom $\Gamma_b$, we use a purely specular reflection  boundary condition. The parameters $\rho_{\text{in}}$, $T_{\text{in}}$, $T_{t}$, $\gamma$, $\nu$ have the same values as in the previous test. We use again~\eqref{eq:IC:2D} as the initial solution.

We note that on the profile of airfoil we cannot use the neighbor points to approximate the tangential derivative $\frac{\partial \hat f}{\partial \hat y}$ in~\eqref{eq:2Dreformulation}. It is because these neighbor points are not on the same straight. Here we approximate the tangential derivative by using the distribution function of interior domain.

In the following tests, we consider only the situations in hydrodynamic regime, {\it i.e.} $\varepsilon=0.05$, for comparing the ones in literature~\cite{bibG, bibK}. A $150\times 100$ mesh is used in domain $\Omega_{\mathbf{x}}$. We use a limit velocity domain $[-8,8]^3$ with mesh size $48\times48\times12$. Two different tests of transonic airflow around this half airfoil are considered: $\text{Mach}_{\text{in}}<1$ and $\text{Mach}_{\text{in}}>1$.

\begin{figure}
 
    \includegraphics[width = 12cm]{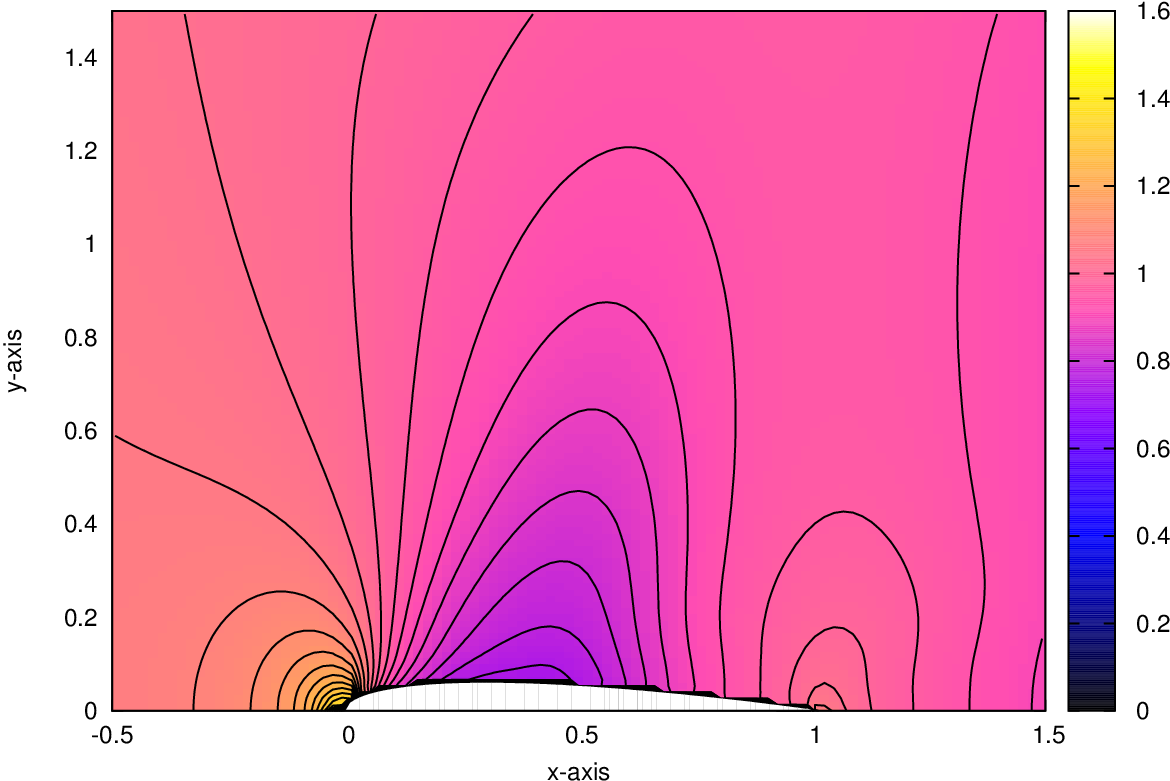} 
 
  \caption{\label{fig:density_wing_M1}Density with $\text{Mach}_{\text{in}}=0.85$}
\end{figure}

\begin{figure}
 
    \includegraphics[width = 12cm]{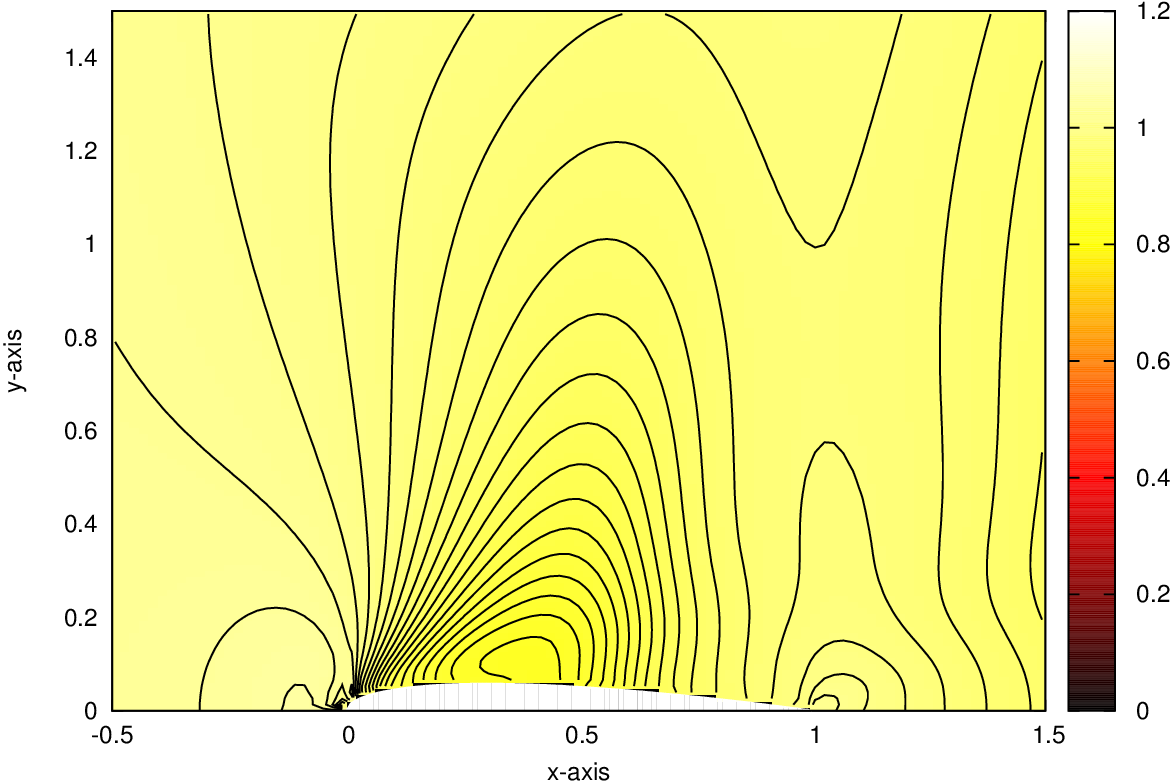} 
 
  \caption{\label{fig:temperature_wing_M1}Temperature with  $\text{Mach}_{\text{in}}=0.85$}
\end{figure}

\begin{figure}
 
    \includegraphics[width = 12cm]{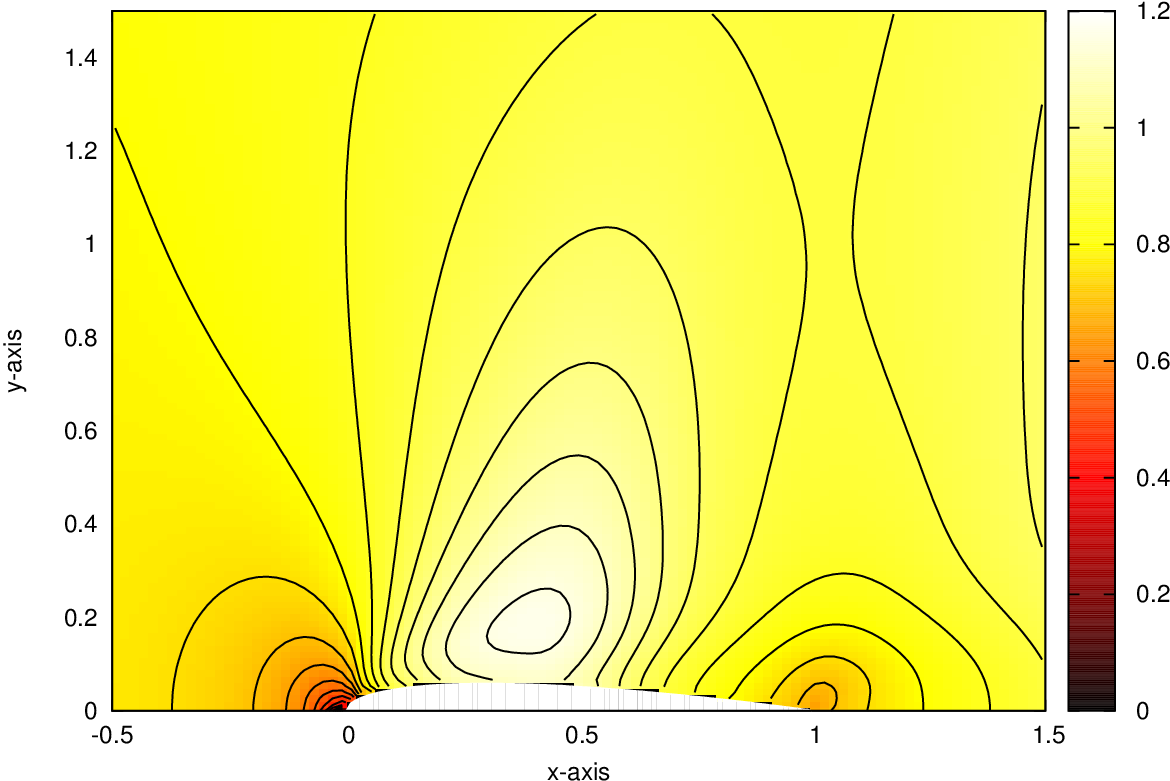} 
 
  \caption{\label{fig:mach_wing_M1}Mach with  $\text{Mach}_{\text{in}}=0.85$}
\end{figure}

\begin{figure}
 
    \includegraphics[width = 12cm]{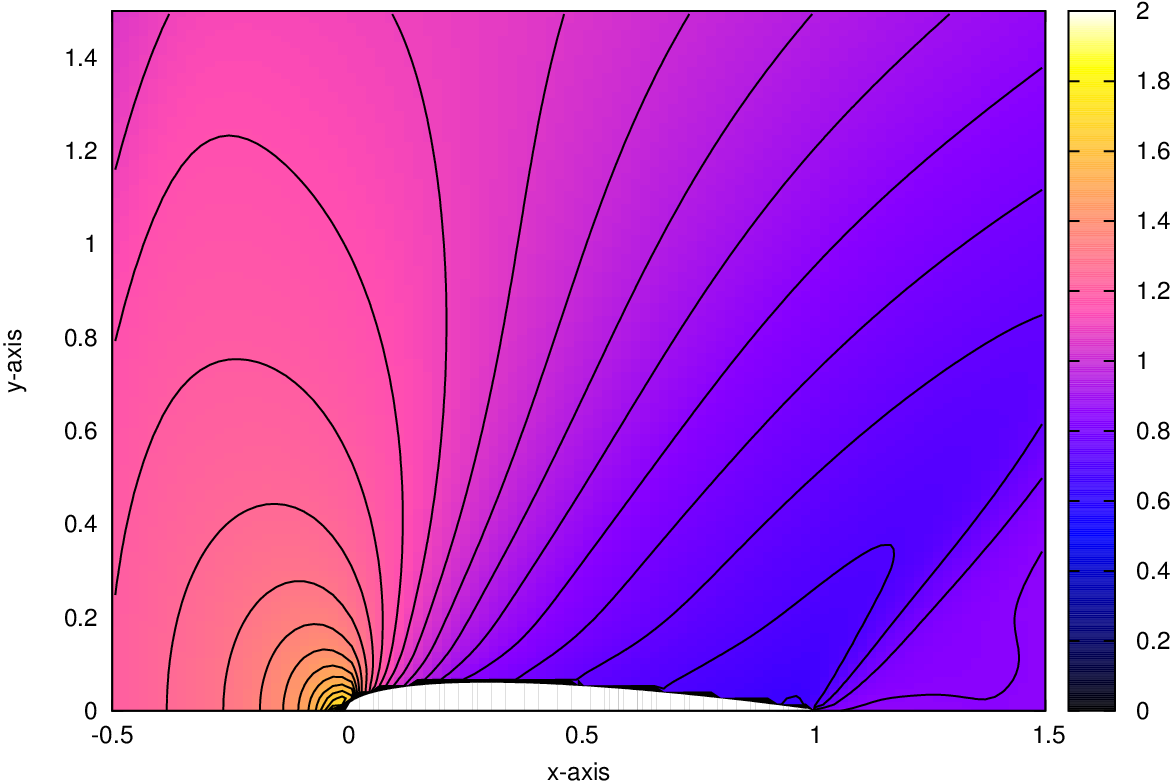} 
 
  \caption{\label{fig:density_wing_M2}Density with $\text{Mach}_{\text{in}}=1.2$}
\end{figure}

\begin{figure}
 
    \includegraphics[width = 12cm]{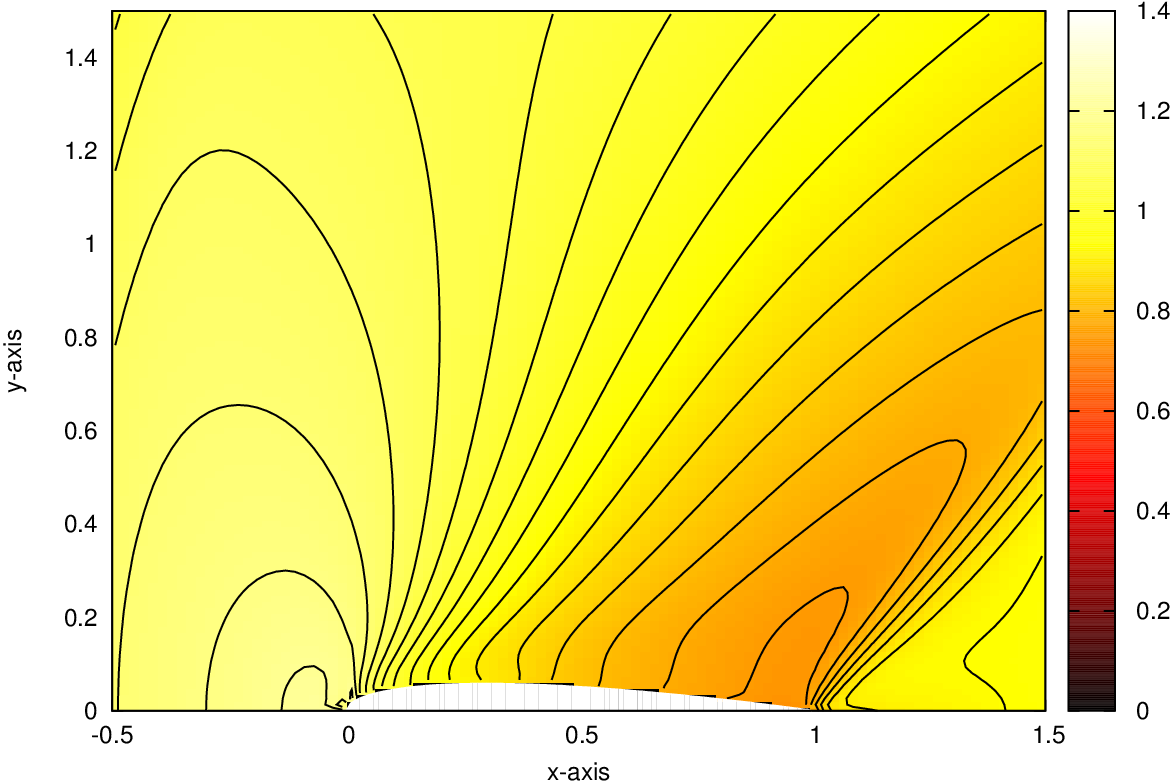} 
 
  \caption{\label{fig:temperature_wing_M2}Temperature with  $\text{Mach}_{\text{in}}=1.2$}
\end{figure}

\begin{figure}
 
    \includegraphics[width = 12cm]{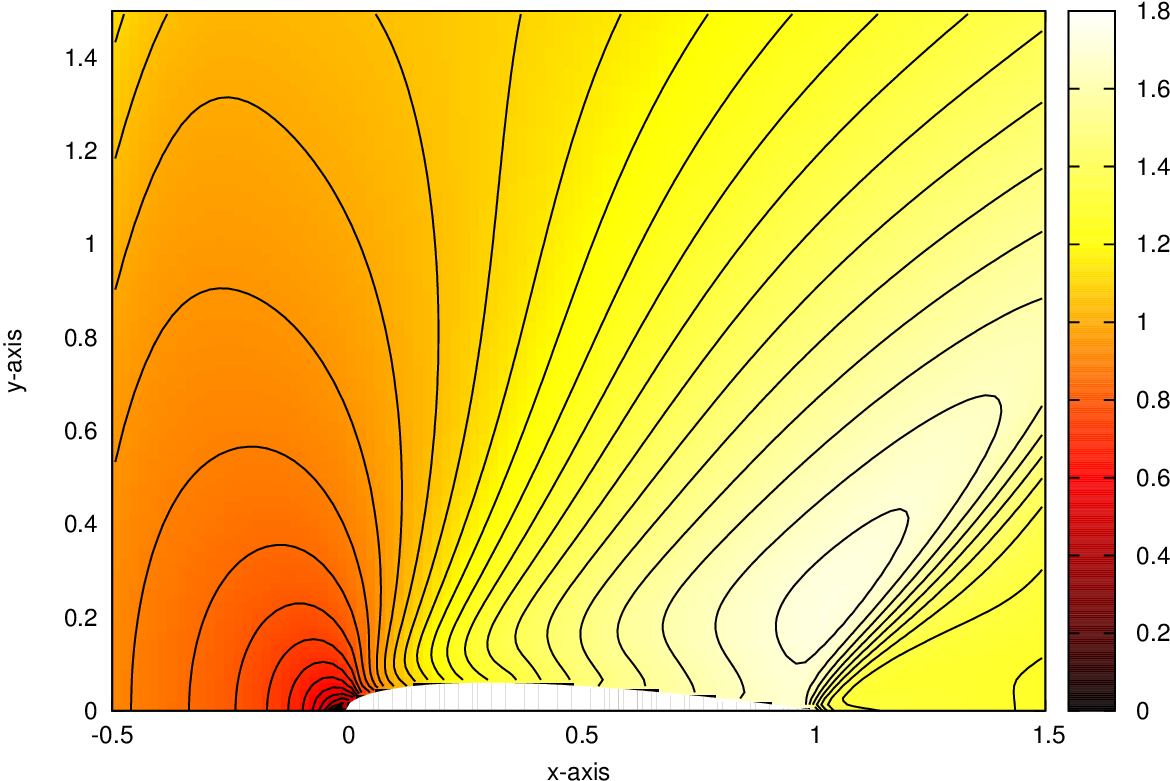} 
 
  \caption{\label{fig:mach_wing_M2}Mach with  $\text{Mach}_{\text{in}}=1.2$}
\end{figure}

We choose first $\text{Mach}_{\text{in}}=0.85$. So we observe in Figures~\ref{fig:density_wing_M1}--\ref{fig:mach_wing_M1} that the flow field around the object includes both sub- ($\text{Mach}<1$) and supersonic ($\text{Mach}>1.2$) parts. The transonic ($0.8<\text{Mach}<1.2$) period begins when first zones of $\text{Mach}>1$ flow appear around the object. Supersonic flow can decelerate back to subsonic before the trailing edge. In the case   $\text{Mach}_{\text{in}}=1.2$ (see in  Figures~\ref{fig:density_wing_M2}--\ref{fig:mach_wing_M2}), the zone of $\text{Mach}>1$ flow increases towards both leading and trailing edges. There is a normal shock created at trailing edge. The flow decelerates over the shock, but remains supersonic. Moreover a normal shock is created ahead of the object, and the only subsonic zone in the flow field is a small area around the  object's leading edge.

\section{Conclusion}
\label{sec:conc}
\setcounter{equation}{0}
In this paper we present an accurate method  based on Cartesian mesh to deal with complex geometry boundary for kinetic models set in a complex geometry.  We desire to reconstruct the distribution function $f$ on some ghost points for  computing transport operator. For this we proceed in three steps: first we extrapolate the distribution function $f$ on  ghost points for outflow. Then we use the boundary conditions to compute the inflow at the boundary. Finally we implement an inverse Lax-Wendroff procedure to give an accurate approximation of $f$ for inflow on the ghost points. A spatially one-dimensional example is given to show that this method has second order accuracy in $L^1$ norm. Moreover some  $1D\times 3D$ and $2D\times 3D$ illustrate that our method can reproduce the similar results as the ones in literature.


\bibliographystyle{plain}


\end{document}